\documentclass[12pt,reqno]{amsart}
\usepackage{amssymb}
\usepackage{amsmath}
\usepackage{amsthm}
\usepackage{eucal}
\usepackage{color}
\usepackage{amsfonts}
\usepackage[pdftex]{graphicx}
\usepackage[T1]{fontenc}

\newcommand{\R}{\mathbb R}

\newcommand{\Z}{\mathbb Z}

\newcommand{\p}{\partial}

\newcommand{\sgn}{\text{sgn}}

\newcommand{\ue}{u^{\epsilon}}
\newcommand{\uel}{u^{\epsilon'}}
\newtheorem{theorem}{Theorem}[section]

\newtheorem{corollary}[theorem]{Corollary}

\newtheorem*{TA1}{Theorem A1}
\newtheorem*{TA2}{Theorem A2}
\newtheorem*{TA3}{Theorem A3}
\newtheorem*{TA4}{Theorem A4}

\numberwithin{equation}{section}
\begin{document}
\title[Regularity of Solutions]{On the regularity of solutions to a class of nonlinear dispersive equations}
\author{Felipe Linares}
\address[F. Linares]{IMPA\\
Instituto Matem\'atica Pura e Aplicada\\
Estrada Dona Castorina 110\\
22460-320, Rio de Janeiro, RJ\\Brazil}
\email{linares@impa.br}

\author{Gustavo Ponce}
\address[G. Ponce]{Department  of Mathematics\\
University of California\\
Santa Barbara, CA 93106\\
USA.}
\email{ponce@math.ucsb.edu}

\author{Derek L. Smith}
\address[D. L. Smith]{Department  of Mathematics\\
University of California\\
Santa Barbara, CA 93106\\
USA.}
\email{ dsl@math.ucsb.edu}

\begin{abstract} We shall study special regularity properties of solutions to some nonlinear dispersive models. The goal is to show how regularity on the initial data is transferred to the solutions. This will depend on the  spaces where regularity is measured.
\end{abstract}

\maketitle

\section{Introduction}

The aim of this work is to study special regularity properties of solutions to the initial value problem (IVP) associated to some nonlinear dispersive
equations of Korteweg-de Vries (KdV) type and related models.

The starting point is a result found by Isaza, Linares and Ponce \cite{ILP-cpde} concerning the solutions of the IVP associated to the 
$k$-generalized KdV equation
\begin{equation}\label{AA}
\begin{cases}
\partial_tu+\partial_x^3u +u^k\partial_xu=0, \hskip5pt \;x, t\in\R, \; k\in\Z^{+},\\
u(x,0)=u_0(x).
\end{cases}
\end{equation}

To state our result we first recall the following local well-posedness for the IVP \eqref{AA} established in \cite{KPV91}:
\begin{TA1}[\cite{KPV91}] 
If $\;u_0\in H^{{3/4}^{+}}(\R)$, then there exist $T\!=\!T(\|u_0\|_{_{{\frac34}^{+}\!, 2}};$ $ k)>0$ and a unique solution of the IVP \eqref{AA}
such that
\begin{equation}\label{notes-2}
\begin{split}
{\rm (i)}\hskip25pt & u\in C([-T,T] : H^{{3/4}^{+}}(\R)),\\
{\rm(ii)}\hskip25pt & \partial_x u\in L^4([-T,T]: L^{\infty}(\R)), \hskip10pt {\text (Strichartz)},\\
{\rm(iii)}\hskip25pt &\underset{x}{\sup}\int_{-T}^{T} |J^r\partial_x u(x,t)|^2\,dt<\infty \text{\hskip10pt for \hskip10pt} r\in[0,{3/4}^{+}],\\
{\rm(iv)}\hskip25pt & \int_{-\infty}^{\infty}\; \sup_{-T\leq t\leq T}|u(x,t)|^2\,dx < \infty,
\end{split}
\end{equation}
with $\;J=(1-\partial_x^2)^{1/2}$. Moreover, the map data-solution, $\,u_0\to \,u(x,t)$ is locally continuous (smooth) 
from $\,H^{3/4+}(\R)$ into the class defined in \eqref{notes-2}.
\end{TA1}

For a detailed discussion on the best  available local and global well-posedness results of the IVP \eqref{AA} we refer to \cite{ILP-cpde} and \cite{LP}.

Now we enunciate the result obtained in \cite{ILP-cpde} regarding propagation of regularities which motivates our study here:

\begin{TA2}[\cite{ILP-cpde}]\label{A2}
If  $u_0\in H^{{3/4}^{+}}(\R)$ and for some $\,l\in \Z^{+},\,\;l\geq 1$ and $x_0\in \R$
\begin{equation}\label{notes-3}
\|\,\partial_x^l u_0\|^2_{L^2((x_0,\infty))}=\int_{x_0}^{\infty}|\partial_x^l u_0(x)|^2dx<\infty,
\end{equation}
then the solution $u=u(x,t)$ of the IVP \eqref{AA} provided by Theorem A1 satisfies  that for any $v>0$ and $\epsilon>0$
\begin{equation}\label{notes-4}
\underset{0\le t\le T}{\sup}\;\int^{\infty}_{x_0+\epsilon -vt } (\partial_x^j u)^2(x,t)\,dx<c,
\end{equation}
for $j=0,1, \dots, l$ with $c = c(l; \|u_0\|_{{3/4}^{+},2};\|\,\partial_x^l u_0\|_{L^2((x_0,\infty))} ; v; \epsilon; T)$.

In particular, for all $t\in (0,T]$, the restriction of $u(\cdot, t)$ to any interval of the form $(a, \infty)$ belongs to $H^l((a,\infty))$.

Moreover, for any $v\geq 0$, $\epsilon>0$ and $R>0$ 
\begin{equation}\label{notes-5}
\int_0^T\int_{x_0+\epsilon -vt}^{x_0+R-vt}  (\partial_x^{l+1} u)^2(x,t)\,dx dt< c,
\end{equation}
with  $c = c(l; \|u_0\|_{_{{3/4}^{+},2}};\|\,\partial_x^l u_0\|_{L^2((x_0,\infty))} ; v; \epsilon; R; T)$.
\end{TA2}

This tells us that the $H^l$-regularity on the right hand side of the data travels forward in time with infinite speed. Notice that since the equation is reversible in time
a gain of regularity in $H^s(\R)$ cannot occur so at $t>0$,  $u(\cdot, t)$ fails to be in $H^j(\R)$ due to its decay at $-\infty$. In
fact, it follows from the proof in  \cite{ILP-cpde} that for any  $\delta>0$ and $t\in (0,T)$ and $j=1,\dots,l$
\begin{equation*}
\int_{-\infty}^{\infty} \frac{1}{\langle x_{-}\rangle^{j+\delta}} (\partial_x^j u)^2(x,t)\,dx \le \frac{c}{t},
\end{equation*}
with $c= c(\|u_0\|_{{3/4}^{+},2}; \|\partial_x^j u_0\|_{L^2((x_0,\infty))}; \,x_0; \,\delta)$.

The result in \cite{ILP-cpde} (Theorem A2) has been extended to the generalized Benjamin-Ono (BO) equation \cite{ILP-bo}
and to the Kadomtsev-Petviashvili II equation \cite{ILP-kp}. Hence, it is natural to ask if this propagation of regularity phenomenon  is intrinsically related to the integrable character of the model or  as in the KdV equation is due to the form of the solution of the associated linear problem.
More precisely, to the structure of its fundamental solution, i.e. the Airy function (see \eqref{1.4b} below). 

Indeed, for the so called $k$-generalized dispersive BO equation,
\begin{equation*}
\p_tu+u^k\p_xu- (-\p_x^2)^{\alpha/2}\p_xu=0, \hskip7pt k\in \Z^{+},\;\;1\le\alpha\le 2,
\end{equation*}
which for $\alpha=1$ corresponds to the $k$-generalized BO equation and $\alpha=2$ to the $k$-generalized KdV equation, one has that the propagation of
regularities (as that presented in Theorem A2) is only known in the cases $\alpha=1$ and $\alpha=2$.
 
Our first result shows that this fact seems to be more general. In particular, it is valid for solutions of the general quasilinear equation of
KdV type, that is,
\begin{equation}\label{i2}
\begin{cases}
\p_tu + a(u,\p_xu,\p_x^2u)\,\p_x^3u + b(u,\p_xu,\p_x^2u)=0,\\
u(x,0)=u_0(x).
\end{cases}
\end{equation}
where the functions $a,b:\R^3\times[0,T]\to \R$ satisfy:
\begin{enumerate}
\item[(H1)]  $\;\;a(\cdot,\cdot,\cdot)$ and  $b(\cdot,\cdot,\cdot)$ are $C^{\infty}$ with all derivatives bounded in $[-M, M]^3$, for
any $M>0$,
\item[(H2)]  given $M>0$, there exists $\kappa>0$ such that 
\begin{equation*}
1/\kappa \le a(x, y, z)\le \kappa \hskip7pt\text{for any}\hskip5pt (x,y,z) \in [-M,M]^3,
\end{equation*}
and
\begin{equation*}
\p_z\,b(x,y,z)\le 0 \hskip5pt\text{for}\hskip5pt (x,y,z) \in [-M,M]^3.
\end{equation*}
\end{enumerate}

To establish the propagation of regularity in solutions of \eqref{i2} of the kind described in \eqref{notes-4} we shall follow the arguments and results obtained by Craig, Kappeler and Strauss in \cite{CKS}.

Assuming the hypotheses (H1) and (H2), the local existence and uniqueness result established in \cite{CKS} affirms:

\begin{TA3}[\cite{CKS}] Let $m\in\Z^{+}$, $m\ge 7$. For any $u_0\in H^m(\R)$, there exist $T=T(\|u_0\|_{7,2})>0$ and a 
unique solution $u=u(x,t)$ of the IVP \eqref{i2} satisfying,
\begin{equation*}
u\in L^{\infty}([0,T] ; H^m(\R)).
\end{equation*}

Moreover, for any $R>0$
\begin{equation*}
\int\limits_0^T\int\limits_{-R}^{R} (\p_x^{m+1} u)^2(x,t)\,dxdt <\infty.
\end{equation*}
 \end{TA3}

For our purpose here we need some (weak) continuous dependence of the solutions upon the data. Hence, we shall first prove the following
\lq\lq refinement\rq\rq \, of Theorem A3.

\begin{theorem}\label{thm1} Let $m\in\Z^{+}$, $m\ge 7$. For any $u_0\in H^m(\R)$ there exist $T=T(\|u_0\|_{7,2})>0$ and a unique
solution $u=u(x,t)$ of the IVP \eqref{i2} such that
\begin{equation}\label{i3}
u\in C([0,T]  : H^{m-\delta}(\R))\cap L^{\infty}([0,T] : H^m(\R)), \hskip5pt \text{for all} \;\;\delta>0,
\end{equation}
with 
\begin{equation}\label{i4}
\p_x^{m+1}u\in L^2([0,T] \times [-R,R]),  \hskip5pt \text{for all} \;\;R>0.
\end{equation}

Moreover, the map data solution $u_0\mapsto u(\cdot, t)$ is locally continuous from $H^m(\R)$ into $C([0,T] : H^{m-\delta}(\R))$ for any
$\delta>0$.
\end{theorem}

Since our objective here is to study propagation of regularities we shall not address the problem of (strong) persistence (i.e. $u_0\in X$, then 
the corresponding solution $u(\cdot,t)$ describes a continuous curve on $X$, $u\in C([0, T] : X)$) and the (strong) continuous dependence
$u_0\mapsto u(\cdot, t)$, (i.e. the map data $\to$ solution from $X$ into $C([0, T] : X)$ is continuous), so that the solutions of \eqref{i2} generate
a continuous flow in $H^m(\R)$, $m\ge 7$.

Our main result concerning the solution of the IVP \eqref{i2} is the following:

\begin{theorem}\label{thm2}
Let $n,m \in\Z^{+}$, $n>m\ge 7$. If $u_0\in H^m(\R)$ and for some $x_0\in\R$
\begin{equation*}
\p_x^j u_0\in L^2((x_0,\infty)) \hskip7pt\text{for}\;\; j=m+1,\dots, n,
\end{equation*}
then the solution of the IVP \eqref{i2} provided by Theorem \ref{thm1} satisfies that for any $\epsilon>0$, $v>0$, and $t\in [0,T)$
\begin{equation}\label{thm2a}
\begin{split}
\int\limits_{x_0+\epsilon-vt}^{\infty}&|\p_x^ju(x,t)|^2\,dx \\
&\le c(\epsilon; v; \|u_0\|_{m,2}; \|\p_x^l u_0\|_{L^2((x_0, \infty))}:l=m+1,\dots,n),
\end{split}
\end{equation}
for $j=m+1, \dots, n$.

Moreover, for any $\epsilon>0$, $v>0$, and $R>0$
\begin{equation}\label{thm2b}
\begin{split}
&\int\limits_{0}^T\int\limits_{x_0+\epsilon-vt}^{x_0+R+vt} |\p_x^{n+1} u(x,t)|^2 \,dxdt \\
&\hskip25pt \le c(\epsilon; v; R;\|u_0\|_{m,2}; \|\p_x^l u_0\|_{L^2((x_0, \infty))}:l=m+1,\dots,n).
\end{split}
\end{equation}
\end{theorem}

Several direct consequences can be deduced from Theorem \ref{thm2} for instance (for further outcomes see \cite{ILP-cpde})

\begin{corollary} Let $u\in C([0,T] : H^m(\R))$, $m\ge 7$, be the solution of the IVP \eqref{i2} provided by Theorem \ref{thm1}. If there exist
$n>m$, $a\in \R$ and $\hat{t}\in (0,T)$ such that 
\begin{equation*}
\p_x^nu(\cdot, \hat{t})\notin L^2((a,\infty)),
\end{equation*}
then for any $t\in (0, \hat{t})$ and any $\beta\in\R$
\begin{equation*}
\p_x^nu(\cdot, t)\notin L^2((\beta,\infty)).
\end{equation*}
\end{corollary}

Theorem \ref{thm2} tells us that the propagation phenomenon described in Theorem A2 still holds in solutions of the quasilinear problem \eqref{i2}. This result and those in KdV, BO, KPII equations seem to indicate that the propagation of regularity phenomena
can be established in systems where Kato smoothing effect (\cite{Kato}) can be proved by integration by parts  directly in the
differential equation.

Since our arguments follow closely those in \cite{CKS} without lost of generality we shall restrict the proofs of Theorem \ref{thm1} and Theorem \ref{thm2} to the case of the model equation
\begin{equation}\label{quasi}
\p_tu + (1+(\p_x^2u)^2) \p_x^3u=0.
\end{equation}

Next we consider the question of the propagation of other type of regularities besides those proved before i.e.  for $u_0\in H^n((x_0, \infty))$ for some 
$x_0\in\R$.

We recall that the next result  can be obtained as a consequence of the argument given by Bona and Saut in \cite{BS}.

\begin{theorem}\label{thm3} Let $k\in \Z^{+}$. There exists 
\begin{equation*}
u_0\in H^1(\R)\cap C^{\infty}(\R)
\end{equation*}
with $\|u_0\|_{1,2}\ll 1$ so that the solution $u(\cdot,t)$ of the IVP \eqref{AA} is
global in time if $k\ge 4$ with
$u\in C(\R: H^1(\R))\cap \dots$ and satisfies
\begin{equation}\label{i10}
\begin{cases}
u(\cdot, t)\in C^1(\R), \;\;t>0,\;\;\;t\notin\Z^{+},\\
u(\cdot, t)\in C^1(\R\backslash \{0\})\backslash C^1(\R), \;\;\; t\in\Z^{+}.
\end{cases}
\end{equation}
\end{theorem}

The argument in \cite{BS} is based in a careful analysis of the asymptotic decay of the Airy function and the well-posedness of the IVP \eqref{AA}
with data $u_0(x)$ in appropriate weighted Sobolev spaces. This argument was simplified (for the case of two points in \eqref{i10}) for the modified KdV equation $k=2$ in \cite{LS} without relying in weighted spaces. Here we shall give a direct proof of Theorem \ref{thm3} which follows the approach in \cite{LS},
i.e.  it does not rely on the analysis of the decay of the Airy function and applies to all the nonlinearities.

Our method has the advantage that it can be extended to $W^{s,p}$-setting. More precisely, we shall show the following:

\begin{theorem}\label{thm4} \hskip10pt
\begin{itemize}
\item[(a)] Fix $k=2, 3, \dots$, let  $p\in (2,\infty)$ and $j\ge 1$, $j\in\Z^{+}$. There exists 
\begin{equation}\label{kdv-lp1}
u_0\in H^{3/4}(\R)\cap W^{j,p}(\R)
\end{equation}
such that the corresponding solution 
\begin{equation*}
u\in C([-T,T]:H^{3/4}(\R))\cap \dots
\end{equation*}
of \eqref{AA} satisfies that there exists $t\in [0,T]$ such that
\begin{equation}\label{kdv-lp2b}
u(\cdot,\pm t)\notin W^{j,p}(\R^{+}).
\end{equation}
\item[(b)] For $k=1$, the same result holds for $j\ge 2$, $j\in \Z^{+}$.
\end{itemize}
\end{theorem} 

\noindent{\bf Remark.} {\em It will follow from our proof that there exists $u_0$ as in \eqref{kdv-lp1} such that \eqref{kdv-lp2b} holds in $\R^{-}$.
Hence, one can conclude that regularities in  $W^{j,p}(\R)$ for $p>2$ do not propagate forward or backward in time to the right or to the left.}

Next  we study the propagation of regularities in solutions to some related dispersive models. First, we consider the IVP associated to the Benjamin-Bona-Mahony BBM equation \cite{BBM}
\begin{equation}\label{ivp-bbm}
\begin{cases}
\p_tu +\p_xu +u\p_xu-\p_x^2\p_tu=0, \hskip7pt x, t\in\R,\\
u(x,0)=u_0(x).
\end{cases}
\end{equation}

The BBM equation was proposed in \cite{BBM} as a model for long surface gravity waves of small amplitude propagating in one dimension.
It was introduced as a "regularized"  version of the KdV equation.   In most cases, the
independent variable $x$ characterizes position in the medium of propagation whilst $t$ is proportional to elapsed time.  The dependent 
variable $u$ may be an amplitude, a pressure, a velocity or other measurable quantity,  depending upon the physical system and the modeling stance taken.

We recall the local well-posedness for the IVP \eqref{ivp-bbm} obtained by Bona and Tzvetkov \cite{BT}.

\begin{TA4}[\cite{BT}] Let $s\ge 0$. For any $u_0\in H^s(\R)$ there exist $T=T(\|u_0\|_{s,2})>0$ and a unique solution $u$ of the IVP \eqref{ivp-bbm}
\begin{equation*}
u\in C([0,T] : H^s(\R)) \equiv  X_T^s.
\end{equation*}

Moreover, the map data-solution $u_0\mapsto u(\cdot, t)$ is locally continuous from $H^s(\R)$ into $X_T^s$.
\end{TA4}

In \cite{BT} it was also shown that Theorem A4 is optimal in an appropriate sense (see \cite{BT} for details).

The following result describes the local propagation of regularities in solutions of the IVP \eqref{ivp-bbm}.

\begin{theorem}\label{thm5}
Let $u_0\in H^s(\R)$, $s\ge 0$. If for some $k\in \Z^{+}\cup \{0\}$, $\theta\in[0,1)$, and $\Omega \subseteq \R$ open
\begin{equation*}
u_0\big|_{\Omega}\in C^{k+\theta},
\end{equation*}
then the corresponding solution $u\in X_T^s$  of the IVP \eqref{ivp-bbm} provided by Theorem A4 satisfies that
\begin{equation*}
u(\cdot, t)\big|_{\Omega} \in C^{k+\theta} \hskip7pt  \text{for all}\;\; t\in [0,T].
\end{equation*}
Moreover, 
\begin{equation*}
u, \p_tu  \in C([0,T]: C^{k+\theta}(\Omega)).
\end{equation*}
\end{theorem}

\underbar{\bf Remarks.}
{\em
\begin{enumerate}
\item Theorem \ref{thm5} tells us that in the time interval $[0,T]$ in the $C^{k+\theta}$ setting no singularities can appear or disappear in the  solution $u(\cdot,t)$.

In particular, one has the following consequence of Theorem \ref{thm5} and its proof:
\begin{corollary}\label{cor-thm5}
Let $u_0\in H^s(\R)$, $s\ge 0$. If for $a<x_0<b$, $k\in \Z^{+}\cup \{0\}$ and $\theta\in [0,1)$
\begin{equation*}
u_0\big|_{(a,x_0)}, \;\;u_0\big|_{(x_0,b)}\in C^{k+\theta} \hskip7pt\text{and}\hskip7pt u_0\big|_{(a,b)}\notin C^{k+\theta},
\end{equation*}
then the corresponding solution $u\in X^s_T$ of the IVP \eqref{ivp-bbm} provided by Theorem A4 satisfies
\begin{equation*}
u(\cdot,t)\big|_{(a,x_0)}, \;\;u(\cdot, t)\big|_{(x_0,b)}\in C^{k+\theta} \hskip7pt\text{and}\hskip7pt u(\cdot, t)\big|_{(a,b)}\notin C^{k+\theta}.
\end{equation*}
\end{corollary}

\item  Theorem A2, Theorem \ref{thm3}, Theorem \ref{thm5}, and Corollary \ref{cor-thm5} show that solutions of the BBM equation  and the KdV equation exhibit a quite different behavior regarding the propagation of regularities.\\
\end{enumerate} 
}
Next we consider the IVP associated to the Degasperis-Procesi (DP) equation (see \cite{DP}):
\begin{equation}\label{ivp-dp}
\begin{cases}
\p_tu-\p_x^2\p_tu+ 4u\p_{x}u=3\p_xu\p_x^2u+u\p_{x}^3u,\hskip7pt x,t\in \R,\\
u(x,0)=u_0(x).
\end{cases}
\end{equation}

The DP model was derived by Degasperis and Procesi as an example of an integrable system similar to the Camassa-Holm (CH) equation 
(\cite{CH})
\begin{equation}\label{c-h}
\p_tu+2\kappa \p_xu-\p_x^2\p_tu+ 3u\p_{x}u=2\p_xu\p_x^2u+u\p_x^3u, \;\;\;\kappa>0.
\end{equation}

Like the  CH equation it possesses  a Lax pair formulation and a bi-Hamiltonian structure leading to an infinite number of conservation laws. 

Also, similar to the CH equation, the DP equation has been shown to exhibit multi-peakons solutions
\begin{equation*}
u(x,t)=\underset{j=1}{\overset{n}{\sum}} \;\alpha_j(t) e^{-|x-x_j(t)|}
\end{equation*}
where $(x_j(t))_{j=1}^n$ satisfies
\begin{equation}\label{peakon}
\begin{cases}
\dfrac{dx_j}{dt}= \underset{k=1}{\overset{n}{\sum}} \alpha_j(t) e^{-|x-x_k(t)|}\\
\\
\dfrac{d\alpha_j}{dt}= 2\alpha_j \,\underset{k=1}{\overset{n}{\sum}} \alpha_k(t) \,{\rm sgn}(x_j-x_k) \,e^{-|x-x_k(t)|}.
\end{cases}
\end{equation}

In \cite{yin} Yin shows that the IVP \eqref{ivp-dp} is locally well-posed in $H^s(\R)$, $s>3/2$.
It shall be remarked that in the Sobolev scale $H^s(\R)$ the result in \cite{yin} is the "optimal" possible for
the "strong" local well-posedness. More precisely, it was established by Himonas, Holliman and
Grayshan \cite{HHG} that for $s<3/2$ solutions of the IVP \eqref{ivp-dp} exhibit an ill-posedness feature due to the so called norm inflection (see \cite{HHG}).

One observes that the multi-peakons solutions barely fail to belong to these spaces since 
$e^{-|x|} \in H^s(\R)$ if and only if $s<3/2$. So we shall construct first a space where uniqueness and existence hold,  where the \lq\lq flow''
(characteristics) is  defined (see \eqref{flow-dp1}  below) and  which contains the multi-peakons.

\begin{theorem}\label{thm6} Given $u_0\in H^{1+\delta}(\R)\cap W^{1,\infty}(\R)$, for some
$\delta>0$. There exist 
$T\!\!=\!T(\|u_0\|_{1+\delta,2}$, $ \|u_0\|_{1,\infty})>0$ and a unique strong solution  (limit of classical solution)
\begin{equation*}
u\in C([0,T]:H^1(\R))\cap L^{\infty}([0,T]:H^{1+\delta}(\R))\cap L^{\infty}([0,T]:W^{1,\infty}(\R)).
\end{equation*}

Moreover, if $u_{n_0}\to u_0$ in $H^{1+\delta}(\R)\cap W^{1,\infty}(\R)$, then the corresponding solutions $u_n, u$ satisfy that $u_n\to u$ in $C([0,T] : H^1(\R))$.

\end{theorem}

With Theorem \ref{thm6} in hand we can show the following result regarding the propagation of regularities in solutions of the DP equation.

\begin{theorem}\label{thm7}  Suppose $u_0\in H^{1+\delta}(\R)\cap W^{1,\infty}(\R)$, for some
$\delta>0$, such that for some open subset  $\Omega\subseteq \R$ 
\begin{equation}\label{thm7-1}
u_0\big|_{\Omega}\in C^1.
\end{equation}
Then the local solution of the IVP \eqref{ivp-dp} $u=u(x,t)$ provided by Theorem \ref{thm6}
satisfies
\begin{equation}\label{thm7-2}
u(\cdot, t)\in C^1(\Omega_t),
\end{equation}
where $\Omega_t=\Phi_t(\Omega)$ and $\Phi_t(x_0)=X(t;x_0)$ is the map given by
\begin{equation}\label{flow-dp1}
\begin{cases}
\dfrac{dX}{dt}= u(X(t),t)\\
\\
X(0)=x_0.
\end{cases}
\end{equation}
\end{theorem}

\noindent\underline{\bf Remarks.}
\begin{enumerate}
{\em
\item Notice that
\begin{equation*}
u\in C([0,T] : H^s(\R))  \hookrightarrow C_b(\R\times [0,T])
\end{equation*}
and
\begin{equation*}
u\in L^{\infty}([0,T] : W^{1,\infty}(\R))
\end{equation*}
which guarantees that the flow in \eqref{flow-dp1} is well defined (see for instance \cite{CL}).
\item Theorem \ref{thm7} describes the propagation of $C^1$ singularities. 
In particular one has the following consequence of Theorem \ref{thm7}.
\begin{corollary}
If $u_0\in H^{1+\delta}(\R)\cap W^{1,\infty}(\R)$, for some $\delta>0$,  (and for some $\Omega\subseteq \R$ open)  $u_0\big|_{\Omega}\in C^1(\Omega\backslash \{\tilde{x}_j\}_{j=1}^{N})$, 
then the corresponding solution $u$  provides by Theorem \ref{thm7}  satisfies for $t\in [0, T]$
\begin{equation*}
u(\cdot, t)\big|_{\Omega_t}\in C^1(\Omega_t\backslash  \{x_j(t)\}_{j=1}^{N})
\end{equation*}
where $\Omega_t$ and $x_j(t)$ are defined by the flow $\Phi_t(x_0)=X(t; x_0)$ as in \eqref{flow-dp1}, $\Omega_t=\Phi(\Omega)$  and 
$x_j(t)=\Phi_t(\tilde{x}_j)$. 
\end{corollary}

\noindent This tell us that all $C^1$-singularities for data in $H^{1+\delta}\cap {\rm Lip}$ propagate with the flow as in the case of multi-peakons  see \eqref{peakon}.

\item Since $m(x,t)= (1-\p_x^2) u(x,t)$ satisfies
\begin{equation}\label{i17}
\p_tm+u\p_xm+3\p_xu\, m=0,
\end{equation}
the propagation at the $C^2$ level and beyond follows directly from the equation  \eqref{i17}. Therefore, the result in Theorem \ref{thm7} extends to $C^k$, $k\in\Z^{+}$, in \eqref{thm7-1} and \eqref{thm7-2}.

\item A result like that  described in Theorem \ref{thm6} (when uniqueness, existence and the
"flow" is defined in a space containing the peakons) is unknown for the CH equation (see
\cite{BC}, \cite{M}, \cite{R}).   Similarly for Theorem \ref{thm7}.
}
\end{enumerate}

Finally,  we consider the 1D version of the Brinkman model \cite{Br}
\begin{equation}\label{brink-orig}
\begin{cases}
\phi\p_t\rho+\p_x(\rho v)= F(t,\rho),\\
(-\tilde{\mu}\p_x^2+\frac{\mu}{\kappa}\big)v=-\p_xP(\rho).
\end{cases}
\end{equation}

This system models fluid flow in certain porous media. It has been useful to treat high porosity systems and a rigidly bounded
porous medium. It also has been used to investigate different convective heat transfer problems in porous media (see \cite{QK} and references
therein).  Here $\rho$ is the fluid's density, $v$ its velocity, $P(\rho)$ is the pressure  and $F$ the external mass flow rate. The physical 
parameters $\mu$, $\kappa$,  $\tilde{\mu}$, and $\phi$ represent  the fluid viscosity, the porous media permeability, the pure fluid viscosity, and the porosity of the media, respectively.  To  simplify we shall assume  (without loss of generality)   
 \begin{equation*}
 \phi=\mu=\kappa=\tilde{\mu}=1\hskip10pt\text{and}\hskip10pt F\equiv 0,
 \end{equation*}
 and that $P(\rho)=\rho^2$

It will be clear from our method of proof that these are not necessary restrictions. Thus the equation in \eqref{brink-orig} becomes
\begin{equation}\label{bk}
\p_x\rho-\p_x(\rho(1-\p_x^2)^{-1}\p_x(\rho^2))=0 \hskip10pt x, t\in\R.
\end{equation}

In \cite{a-iorio}  Arbieto and Iorio established the local well-posedness of the associated IVP to \eqref{bk}
in $H^s(\R)$, $s>3/2$.

Our first outcome in this regard shows that the result in \cite{a-iorio} can be improved to $H^s(\R)$, 
$s\ge 1$.

\begin{theorem}\label{thm8} Let $\rho_0\!\in\! H^s(\R)$, $s\ge1$. \!There exist 
$T=T(\|u_0\|_{1,2})=c\,\|u_0\|_{1,2}^{-2}$ and a unique solution 
$\rho=\rho(x,t)$ such that
\begin{equation*}
\rho\in C([0,T] : H^s(\R))\cap C^1((0,T) : H^{s-1}(\R)).
\end{equation*}
Moreover, the map data solution is locally continuous from $H^s$ to  $C([0,T] : H^s(\R))$.
\end{theorem}

With this local well-posedness  theory we establish the following result concerning the propagation of regularities in solutions of the IVP associated to \eqref{bk}.

\begin{theorem}\label{thm9} Suppose $\rho_0\in H^s(\R)$, $s\ge 1$, such that
\begin{equation*}
\rho_0\big|_{\Omega}\in C^k(\Omega) \hskip7pt \text{for some}\;\; \Omega \subseteq \R \hskip5pt\text{open}.
\end{equation*}
Then the corresponding solution $\rho=\rho(x,t)$ of the IVP associated to \eqref{bk} provided by Theorem \ref{thm8} satisfies that
\begin{equation*}
\rho(\cdot, t)\big|_{\Omega_t} \in C^k(\Omega_t),
\end{equation*}
where $\Omega_t=\Phi_t(\Omega)$ and $\Phi_t(x_0)=X(t;x_0)$ is the map defined by
\begin{equation*}
\begin{cases}
\dfrac{dX}{dt}= u(X(t),t),\\
\\
X(0)=x_0.
\end{cases}
\end{equation*}
\end{theorem}

The rest of this work is organized as follows: Theorem \ref{thm1} and Theorem \ref{thm2}
will be proven in Section 2. Section 3 contains the proof of Theorem \ref{thm3} and
Theorem \ref{thm4}. Section 4 is concerned with Theorem \ref{thm5}, Section 5 with
Theorem \ref{thm6} and Theorem \ref{thm7} and Section 6 with Theorem \ref{thm8}
and Theorem \ref{thm9}.

\section{Quasilinear KdV type equations}

For simplicity in the presentation and without loss of generality we will consider the following IVP associated to a quasilinear KdV type
satisfying the hypotheses (H1)-(H2) in \eqref{i2}

\begin{equation}\label{ql-kdv1}
\begin{cases}
\p_tu +\big(1+(\p_x^2u)^2\big)\p_x^3u=0,\hskip5pt x\in\R,\;t>0,\\
u(x,0)=u_0(x)\in H^m(\R), \;\;\; m\ge 7.
\end{cases}
\end{equation}
\vspace{3mm}

\noindent\underline{\bf Proof of Theorem \ref{thm1}.}
\vspace{3mm}

For $\epsilon\in (0,1)$ consider the parabolic problem
\begin{equation}\label{ql-kdv2}
\begin{cases}
\p_tu +\big(1+(\p_x^2u)^2\big)\p_x^3u=-\epsilon\p_x^4u,\\
u(x,0)=u_0(x)\in H^7(\R).
\end{cases}
\end{equation}

Denoting by
\begin{equation*}
K_t\ast f= e^{-\epsilon t\p_x^4}f=\big(e^{-\epsilon t(2\pi\xi)^4}\widehat{f} \,\big)^{\vee}
\end{equation*}
using that
\begin{equation*}
\|\p_x^l(K_t\ast f)\|_2\le \frac{c_l}{(\epsilon t)^{l/4}}\,\|f\|_2,
\end{equation*}
and writing the solution $u^{\epsilon}$ of \eqref{ql-kdv2} in its equivalent integral equation version
\begin{equation}\label{ql-kdv3}
u^{\epsilon}(t)=K_t\ast u_0+\int\limits_0^t K_{t-t'}\ast \big(1+(\p_x^2u^{\epsilon})^2\big)\p_x^3u^{\epsilon}(t')\,dt',
\end{equation}
one sees that there exist 
\begin{equation*}
T_{\epsilon}=T_{\epsilon}(\epsilon; \|u_0\|_{m,2})>0\hskip5pt\text{with}\hskip5pt T_{\epsilon}\sim {\rm O}(\epsilon^{3/4})
\end{equation*}
and a unique $u^{\epsilon}$ solution of \eqref{ql-kdv2} such that
\begin{equation}\label{ql-kdv4}
u^{\epsilon}\in C([0, T_{\epsilon}] : H^m(\R))\cap  C^{\infty}((0, T_{\epsilon}) : H^{\infty}(\R)).
\end{equation}

Next we apply the argument in \cite{CKS} to the IVP \eqref{ql-kdv2} to obtain {\it  a priori} estimate which allows us to extend the local solution $u^{\epsilon}(\cdot)$ in a time
interval $[0,T]$ with $T=T(\|u_0\|_{m,2})>0$ independent of $\epsilon$ in the  class  described in \eqref{ql-kdv4}. The only difference with the argument provided in \cite{CKS} is the term on the 
RHS of the equation in \eqref{ql-kdv2}. This can be easily handled using that
\begin{equation}\label{ql-kdv5}
\begin{split}
-\epsilon \int \p_x^{4+l}\ue\p_x^l \ue \chi_l(x,t)\,dx &= -\epsilon\int \p_x^{2+l}\ue \p_x^{2+l}\ue \chi_l(x,t)\,dx\\
&\;\;\;\;-2\epsilon\int \p_x^{2+l}\ue\p_x^{1+l}\ue \p_x\chi_l(x,t)\,dx\\
&\;\;\;\; -\epsilon\int \p_x^{2+l}\ue \p_x^{l}\ue \p_x^2\chi_l(x,t)\,dx\\
&= -\epsilon\int (\p_x^{2+l}\ue)^2 \chi_l(x,t)\,dx\\
&\;\;\;\; -2\epsilon\int (\p_x^{2+l}\ue)(\p_x^l\ue)\p_x^2\chi_l(x,t)\,dt\\
&\;\;\;\;+\frac{\epsilon}2\int (\p_x^l\ue)(\p_x^l\ue)\p_x^4\chi_l(x,t)\,dt\\
&= E_1+E_2+E_3,
\end{split}
\end{equation}
with $\chi_l(u)$ such that
\begin{equation}\label{ql-kdv6}
\frac32\p_x [(1+(\p_x^2u)^2) \chi_l]-(l+1)\p_x(1+(\p_x^2u)^2)\chi_l=0,\;\;\;l\in \Z^{+},
\end{equation}
i.e. 
\begin{equation*}
\chi_l= (1+(\p_x^2\ue)^2)^{-c_l} , \hskip15pt  c_l=\frac23(l+1)-1.
\end{equation*}

Thus one sees that $E_1$ has the appropriate sign,
\begin{equation*}
\begin{split}
E_2 &\le \frac{\epsilon}{2} \int (\p_x^{2+l}\ue)^2 \chi_l(x,t)\,dx + 2\epsilon\int (\p_x^l\ue)^2 \frac{(\p_x^2\chi_{l})^2}{\chi_l}(x,t)\,dx\\
&= E_{2,1} + E_{2,2},
\end{split}
\end{equation*}
where $E_{2,1}$ is absorbed by $E_1$.  Terms  $E_{2,2}$ and $E_3$  are at the level of the estimate involving derivatives of order $l$.

Hence, we can conclude that there exist $T=T(\|u_0\|_{7,2})>0$ and a unique solution 
\begin{equation*}
\ue\in C([0,T]:H^7(\R))\cap C((0,T):H^{\infty}(\R))
\end{equation*}
of \eqref{ql-kdv2} with
\begin{equation}\label{ql-kdv7}
\underset{0<\epsilon<1}{\sup}\; \underset{0\le t\le T}{\sup} \|\ue(t)\|_{7,2}\le M_7= M_7(\|u_0\|_{7,2}).
\end{equation}

Moreover, if $u_0\in H^m(\R)$ with $m\ge 7$, $m\in\Z^{+}$, then
\begin{equation*}
\ue \in C([0,T]:H^m(\R))\cap C((0,T):H^{\infty}(\R))
\end{equation*}
with
\begin{equation*}
\underset{0<\epsilon<1}{\sup}\; \underset{0\le t\le T}{\sup} \|\ue(t)\|_{m,2}\le M_m= M_m(\|u_0\|_{m,2}).
\end{equation*}

Next we want to show that $(\ue)_{\epsilon>0}$ converges in $L^{\infty}([0,T] :H^3(\R))$.

Thus we consider
\begin{equation*}
\p_t\ue +(1+(\p_x^2\ue)^2)\p_x^3\ue= -\epsilon \p_x^4\ue
\end{equation*}
and 
\begin{equation*}
\begin{split}
\p_t\p_x^3\ue +(1+(\p_x^2\ue)^2)\p_x^6\ue+4\p_x(1+(\p_x^2\ue)^2)\p_x^5\ue&\\
+6\p_x^2\ue\p_x^4\ue\p_x^4\ue+12\p_x^3\ue\p_x^3\ue\p_x^4\ue&= -\epsilon \p_x^7\ue
\end{split}
\end{equation*}

Similarly for $\uel$ with $0<\epsilon'<\epsilon<1$. Now we write the equations for $w=w^{\epsilon,\epsilon'}=\ue-\uel$
and $\p_x^3 w$, i.e.
\begin{equation}\label{ql-kdv8}
\begin{split}
\p_tw&+(1+(\p_x^2\ue)^2)\p_x^3w+\p_x^3\uel(\p_x^2\ue+\p_x^2\uel)\p_x^2w\\
&=-\epsilon'\p_x^3w-(\epsilon-\epsilon')\p_x^4\uel
\end{split}
\end{equation}
and
\begin{equation}\label{ql-kdv9}
\begin{split}
\p_t\p_x^3w&+(1+(\p_x^2\ue)^2)\p_x^6w+(\p_x^6\uel)(\p_x^2\ue+\p_x^2\uel)\,\p_x^2w\\
&+8\p_x^2\ue\p_x^3\ue\p_x^5w+4\p_x^4\uel\,\p_x\big[ (\p_x^2\ue+\p_x^2\uel)\p_x^2w\big]\\
&+6\p_x^2\ue(\p_x^4\ue+\p_x^4\uel)\p_x^4w +6(\p_x^4\uel)^2\p_x^2w\\
&+12(\p_x^3\ue)^2\p_x^4w+12(\p_x^4\uel)(\p_x^3\ue+\p_x^3\uel)\p_x^3w\\
&=-\epsilon' \p_x^7w-(\epsilon-\epsilon')\p_x^7\uel.
\end{split}
\end{equation}

We multiply \eqref{ql-kdv8} by $w$ and integrate the result to get
\begin{equation*}
\begin{split}
&\frac12\frac{d}{dt}\int w^2+\frac32\int \p_x((\p_x^2\ue)^2)(\p_x w)^2-\frac12\int \p_x^3((\p_x^2\ue)^2)w^2\\
&-\int (\p_x^3\ue)(\p_x^2\ue+\p_x^2\uel)(\p_xw)^2+\frac12\int \p_x^2\big[\p_x^3\uel(\p_x^2\ue+\p_x^2\uel)\big]w^2\\
&=-(\epsilon-\epsilon')\int \p_x\uel w
\end{split}
\end{equation*}
and multiply \eqref{ql-kdv9} by $\p_x^3w\,\chi$ and integrate the result to get (after some integration by parts)
\begin{equation*}
\begin{split}
&\frac12\frac{d}{dt}\int (\p_x^3w)^2\chi-\int(\p_x^3w)^2\p_x\chi\\
&+\int \Big\{\frac32 \p_x\big[(1+(\p_x^2\ue)^2)\chi\big]-(3+1)\p_x(1+(\p_x^2\ue)^2)\,\chi\Big\} (\p_x^4w)^2\\
&-\frac12\int\p_x^3(1+(\p_x^2\ue)^2\chi)(\p_x^3w)^2+4\int \p_x^2(\p_x^2\ue\p_x^3\ue \chi)(\p_x^3w)^2\\
&-3\int\p_x[\p_x\ue(\p_x^4\ue+\p_x^4\uel)\chi] (\p_x^3w)^2 -6\int \p_x[(\p_x^3\ue)^2\chi] (\p_x^3w)^2\\
&+\int \p_x^6\uel(\p_x^2\ue +\p_x^2\uel)\p_x^2w\p_x^3w+ 4\int\p_x^5\uel \p_x\big((\p_x^2\ue+\p_x^2\uel)\p_x^2w\big)\p_x^3w\chi\\
&+6\int (\p_x^4\uel)^2\p_x^2w\p_x^3w \chi+  12\int (\p_x^3\ue+\p_x^3\uel)\p_x^4\uel (\p_x^3w)^2\chi\\
&=-\epsilon' \int \p_x^7w\p_x^3w \chi -(\epsilon-\epsilon') \int\p_x^7\uel \p_x^3w \chi.
\end{split}
\end{equation*}
As in \cite{CKS} one chooses $\chi=\chi_{3}$ (see \eqref{ql-kdv6}) such that
\begin{equation*}
\frac32 \p_x[(1+(\p_x^2\ue)^2)\chi]-(3+1)\p_x\big(1+(\p_x^2\ue)^2\big)\chi\equiv 0,
\end{equation*}
i.e. 
\begin{equation*}
c_l=\frac23(l+1)-1>0
\end{equation*}
and
\begin{equation*}
\chi_3=(1+(\p_x^2\ue)^2)^{-c_3}\sim 1,
\end{equation*}
so by hypothesis and previous results in \eqref{ql-kdv7} it follows that
\begin{equation*}
\underset{[0,T]}{\sup} \|(\ue-\uel)(t)\|_{3,2}\le c\big(\|u_0\|_{7,2}\big)\Big[ \|u_0^{\epsilon}-u_0^{\epsilon'}\|_{3,2}+(\epsilon+\epsilon')\Big].
\end{equation*}

Hence
\begin{equation*}
(\ue)_{\epsilon}\subseteq C([0,T] : H^7(\R))
\end{equation*}
(uniformly bounded) converges as $\epsilon\downarrow 0$ in $L^{\infty}([0,T] : H^3(\R))$ to $u\in C([0,T]: H^3(\R))$.

Moreover, from \eqref{ql-kdv7} it follows that for any $\delta>0$
\begin{equation*}
\ue \to u \hskip7pt \text{in}\hskip7pt C([0,T] : H^{7-\delta}(\R))
\end{equation*}
and
\begin{equation*}
u\in L^{\infty}([0,T] : H^7(\R))
\end{equation*}
is a solution of the IVP \eqref{ql-kdv1}. But by the uniqueness established in \cite{CKS} this solution is unique  so it agrees with that provided in Theorem A3.

In particular, for any $R>0$
\begin{equation}\label{A0b}
\int\limits_0^T\int\limits_{-R}^{R} |\p_x^8 u(x,t)|^2\,dxdt \le c = c(T; R; \|u_0\|_{7,2}).
\end{equation}

Moreover, the previous argument applied to two solutions 
\begin{equation*}
u, v\in C([0,T] : H^{7-\delta}(\R))\cap L^{\infty}([0,T] : H^7(\R))
\end{equation*}
with data $u_0, v_0$ respectively, shows that
\begin{equation*}
\underset{[0,T]}{\sup} \|(u-v)(t)\|_{3,2} \le c\big( M(\|u_0\|_{7,2},\|v_0\|_{7,2})\big) \|u_0-v_0\|_{3,2}.
\end{equation*}

Hence, the map data solution from $H^7(\R)$ into $C([0,T] : H^{7-\delta}(\R))$ is locally continuous for any $\delta>0$.

\vspace{5mm}
\noindent\underline{\bf Proof of Theorem \ref{thm2}.}

\vspace{3mm}

As in the previous section we will consider the following IVP
\begin{equation}\label{ql-kdv1b}
\begin{cases}
\p_tu +\big(1+(\p_x^2u)^2\big)\p_x^3u=0,\hskip5pt x\in\R,\;t>0,\\
u(x,0)=u_0(x)\in H^m(\R), \;\;\; m\ge 7.
\end{cases}
\end{equation}

Without loss of generality we shall assume that $m=7$ and that $x_0=0$. Thus we have
that the solution $u(\cdot)$ of the IVP \eqref{ql-kdv1b} provided by Theorem \ref{thm1} satisfies
\begin{equation}\label{A1}
u\in C([0,T] : H^{7-\delta}(\R))\cap L^{\infty}([0,T] :  H^7(\R)), \hskip7pt\text{for all}\;\; \delta>0
\end{equation}
and \eqref{A0b}.

As in \cite{ILP-cpde} we define the family of cut off functions: for $\epsilon>0$ and $b\ge 5\epsilon$, 
let $\chi_{_{\epsilon, b}}\in C^{\infty}(\R)$ such that
\begin{equation}\label{A3}
\chi_{_{\epsilon, b}}=
\begin{cases}
0, \hskip5pt x\le\epsilon,\\
1, \hskip5pt x\ge b,
\end{cases}
\end{equation}

\begin{equation}\label{A4}
{\rm supp}\,\chi_{_{\epsilon, b}} \subseteq [\epsilon, \infty), \hskip5pt
{\rm supp}\,\chi_{_{\epsilon, b}}' \subseteq [\epsilon, b],
\end{equation}
\begin{equation}\label{A5}
\chi_{_{\epsilon, b}}'(x)\ge 0 \hskip7pt\text{with}\hskip7pt \chi_{_{\epsilon, b}}'(x)\ge \frac{1}{b-3\epsilon}
1_{[3\epsilon, \,b-2\epsilon]}(x).
\end{equation}

Thus,
\begin{equation}\label{A6}
\chi_{_{\epsilon/3, b+\epsilon}}'(x)\ge c_j\,|\chi_{_{\epsilon, b}}^{(j)}(x)|, \hskip10pt x\in\R, \;j=1, 2, 3.
\end{equation}

Next we follow the argument in \cite{CKS}. Thus we formally apply $\p_x^j$, $j=8, \dots, m$ to the
equation in \eqref{ql-kdv1b} and multiply the result by
\begin{equation}\label{A7}
\p_x^ju \, \p_x^j\psi_j=\p_x^ju \,\psi_{j,v,\epsilon,b},
\end{equation}
with $ \psi_{j,v,\epsilon,b}(x,t)$ to be determine below, to get (after some integration by parts)
\begin{equation}\label{A8}
\begin{split}
&\frac12 \frac{d}{dt} \int (\p_x^ju)^2 \psi_j(x,t)\,dx -\frac12\int (\p_x^ju)^2 \p_t\psi_j(x,t)\,dx\\
&+\int (\p_x^{j+1}u)^2\Big\{ \frac32 \p_x\big((1+(\p_x^2u)^2\big)\psi_j\big)
-j\big[ \p_x(1+(\p_x^2u)^2)\psi_j\big]\\
&\;\;\;\;-\big[\p_x(1+(\p_x^2u)^2)\psi_j\big]\Big\}\,dx\\
&+\int Q_j\big ( (\p_x^lu)_{|l|\le j} ; (\p_x^r\psi_j)_{|r|\le3}\big)(x,t)\,dx\\
&\equiv \frac12  \frac{d}{dt} \int (\p_x^ju)^2 \psi_j(x,t)\,dx +E_{1,j}+E_{2,j}+E_{3,j}.
\end{split}
\end{equation}
where $Q_j(\cdot, \cdot)$ is a polynomial in its variables, linear in the components $(\psi_j, \p_x\psi_j,
\p_x^2\psi_j, \p_x^3\psi_j)$, and at most quadratic in the highest derivatives of $u$, i.e. $(\p_x^ju)$, involving at most $7+2j$ derivatives of $u$.

First we consider $E_{2,j}$ which determines the choices of $\psi_j=\psi_{j,v,\epsilon,b}$. As in
\cite{CKS} we choose  $\psi_j=\psi_{j,v,\epsilon,b}$ such that for $v>0$
\begin{equation}\label{A9}
\frac32 \p_x\big((1+(\p_x^2u)^2\big)\psi_j\big)-(j+1) \p_x\big(1+(\p_x^2u)^2\big)\psi_j
= \chi_{_{\epsilon, b}}' (x+vt)
\end{equation}
with
\begin{equation}\label{A10}
\psi_j(x,t) \to 0 \hskip10pt\text{as}\hskip10pt x\downarrow -\infty,
\end{equation}
i.e.
\begin{equation}\label{A11}
\begin{split}
(1+(\p_x^2u)^2)\p_x\psi_j -&\Big(1-\frac23(j+1)\Big) \p_x\big(1+(\p_x^2u)^2\big)\psi_j\\
&= \frac23\chi_{_{\epsilon, b}}' (x+vt).
\end{split}
\end{equation}

Hence if
\begin{equation}\label{A12}
d_j \equiv 1-\frac23(j+1),
\end{equation}
then
\begin{equation*}\label{A13a} 
\begin{split}
&\psi_j=\psi_{j,v,\epsilon,b}\\
&\equiv \frac23 (1+(\p_x^2u)^2)^{-d_j}(x,t)\! \int\limits_{-\infty}^x \!(1+(\p_x^2u)^2)^{d_j-1}(s,t) \chi_{_{\epsilon,b}}'(s+vt)\,ds.
\end{split}
\end{equation*}

We  observe that
\begin{equation}\label{A13b}
\psi_j(x,t)=\psi_{j,v,\epsilon,b}(x,t)\ge 0
\end{equation}
with
\begin{equation}\label{A13c}
{\rm supp}\, \psi_{j,v,\epsilon,b}(\cdot, t)\subseteq [\epsilon-vt,\infty), \hskip10pt t\in[0,T]
\end{equation}
and 
\begin{equation}\label{A13d}
\psi_{j,v,\epsilon,b}(\cdot, t)\in L^{\infty}(\R)
\end{equation}
with
\begin{equation}\label{A13e}
\|\psi_{j,v,\epsilon,b}(\cdot, t)\|_{\infty} \le c=c(v,\epsilon, b;  \|u_0\|_{7,2}).
\end{equation}

With this choice of $\psi_j(\cdot)=\psi_{j,v,\epsilon,b}(\cdot)$, $E_{2,j}$  becomes 
\begin{equation}\label{A14}
E_{2,j}=\int (\p_x^{j+1}u)^2(x,t)  \chi_{_{\epsilon,b}}'(x+vt)\, dx.
\end{equation}

Also one sees (using \eqref{A1}-\eqref{A0b} and (H2) ) that for any $T>0$, $v>0$, $j\in\Z^{+}$, $j\ge 8$, there exists
$c= c(T; v; j; k)\in (1,\infty)$ such that for any $(x,t)\in \R\times[0,T]$
\begin{equation}\label{A15}
c^{-1}  \chi_{_{\epsilon,b}}(x+vt) \le \psi_{j,v,\epsilon,b}(x,t) \le c\, \chi_{_{\epsilon,b}}(x+vt).
\end{equation}
Moreover (with $c$ as above)
\begin{equation}\label{A16}
\begin{split}
|\p_x\psi_{j,v,\epsilon, b}(x,t)| &\le \big(\chi_{_{\epsilon,b}}(x+vt) + \chi_{_{\epsilon,b}}'(x+vt)\big)\\
&\le  c \big( \psi_{j,v,\epsilon, b}(x,t)+\chi_{_{\epsilon,b}}'(x+vt)\big),
\end{split}
\end{equation}
\begin{equation}\label{A17}
\begin{split}
|\p_x^2\psi_{j,v,\epsilon, b}(x,t)| &\le \big( (\chi_{_{\epsilon,b}}+\chi_{_{\epsilon,b}}')(x+vt)\big)\\
& \;\;\;\;+ |\chi_{_{\epsilon,b}}''(x+vt)|\big),
\end{split}
\end{equation}
and
\begin{equation}\label{A18}
\begin{split}
|\p_x^3\psi_{j,v,\epsilon, b}(x,t)| &\le \big(\big(\chi_{_{\epsilon,b}}+\chi_{_{\epsilon,b}}')(x+vt)\\
&\;\;\;\; + |\chi_{_{\epsilon,b}}''(x+vt)|+ |\chi_{_{\epsilon,b}}^{(3)}(x+vt)|\big),
\end{split}
\end{equation}

Therefore, combining \eqref{A15} and \eqref{A6} one has that 
\begin{equation}\label{A19}
|\p_x^r \psi_{j,v,\epsilon,b}(x,t)|\le c\big(\psi_{j,v,\epsilon,b}(x,t)+\chi_{_{\epsilon,b}}'(x+vt)+\chi_{_{\epsilon/3,b+\epsilon}}(x+vt),
\end{equation}
$r=2, 3$, for any $(x,t)\in\R\times[0,T]$.

Also using the equation, \eqref{A1}, and \eqref{A15}  it follows that
\begin{equation}\label{A20}
\begin{split}
|\p_t\psi_{j,v,\epsilon,b}(x,t)| &\le c \big(\chi_{_{\epsilon,b}}(x+vt)+\chi_{_{\epsilon,b}}'(x+vt)\big)\\
&\le c \big(  \psi_{j,v,\epsilon,b}(x,t)+ \chi_{_{\epsilon,b}}'(x+vt)\big).
\end{split}
\end{equation}

We now turn to the estimate of $E_{1,j}$ in \eqref{A8}. First, we consider the case $j=8$. In this case using \eqref{A20}
we have that 
\begin{equation}\label{A20b}
\begin{split}
|E_{1,8}| &\le \frac{c}{2} \int (\p_x^8u)^2 \psi_{j,v,\epsilon,b} (x,t)\, dx\\
&\;\;\; + c\,\big |\int (\p_x^8u)^2  \chi_{_{\epsilon,b}}'(x+vt)\, dx\big| = E_{1,8,1}+E_{1,8,2}.
\end{split}
\end{equation}

We notice that $E_{1,8,1}$ is a multiple of the term we are estimating, so it will be part of the
Gronwall's inequality. For the term $E_{1,8,2}$ we observe that given $T>0$ and $v>0$, there
exists $R>0$ such that
\begin{equation}\label{A21}
 \chi_{_{\epsilon,b}}'(x+vt)\le 1_{(-R,R)} (x)
 \end{equation}
 for all $(x,t)\in \R\times [0,T]$. Therefore, the bound of $E_{1,8,2}$ follows from \eqref{A0b}
 after integrating in time in the estimate. So at the level $j=8$ it only remains to
 consider $E_{3,8}$ in \eqref{A8}.
 
 Using the structure of $E_{3,8}$ commented after \eqref{A8} and the bounds in 
 \eqref{A15}-\eqref{A19} it follows that  in this case $j=8$, the term $E_{3,8}$ can be written as 
 a sum of terms of the form
 \begin{equation}\label{A22}
 \int  P_{3,8,1} \big((\p_x^lu)_{| l |\le 4} \big)\, (\p_x^8u)^2\, \psi_{8,v,\epsilon,b}(x,t) \,dx \;\;\; (\equiv E_{3,8,1}),
 \end{equation}
 \begin{equation}\label{A23}
 \int P_{3,8,2,r}\big((\p_x^l u)_{|l|\le 4-r}\big) \,(\p_x^8u)^2\,  \chi_{_{\epsilon,b}}^{(r)}(x+vt) \,dx \;\;\; (\equiv E_{3,8,2})
 \end{equation}
with $r=1, 2, 3$, or terms involving lower order derivatives $l=0,\dots,7$,
\begin{equation}\label{A24}
 \int P_{3,8,3}\big((\p_x^l u)_{|l|\le 7}\big)\, \chi_{_{\epsilon,b}}^{(r)}(x+vt) \,dx \;\;\; (\equiv E_{3,8,3}),
 \end{equation}
\begin{equation}\label{A25}
 \int P_{3,8,4}\big((\p_x^l u)_{|l|\le 7}\big)\,  \psi_{8,v,\epsilon,b}(x, t) \,dx \;\;\; (\equiv E_{3,8,4}),
 \end{equation}
 where  
 \begin{itemize}
 \item[ ] $ P_{3,8,1} (\cdot)$ and  $P_{3,8,2,r} (\cdot)$ are polynomials of degree at most 2, 
 \end{itemize}
 and
  \begin{itemize}
 \item[ ]  $P_{3,8,3} (\cdot)$ and  $P_{3,8,4} (\cdot)$ are polynomials of degree at most 4.
  \end{itemize}
  
 In this case, $j=8$, the terms in $E_{3,8,3}$ and $E_{3,8,4}$ are bounded since one has \eqref{A1} and the
 fact that
 \begin{equation*}
 \underset{r=3}{\overset{3}{\sum}}\; \| \chi_{_{\epsilon,b}}^{(r)}\|_{\infty}
 +\|\psi_{8,v,\epsilon,b}\|_{L^{\infty}(\R\times[0,T])}\le c,
 \end{equation*}
 with $c = c(v,\epsilon, b, \underset{[0,T]}{\sup} \|u(t)\|_{7,2})$.
 
 The term $E_{3,8,2}$ can be estimated using \eqref{A0b} as in \eqref{A21}.
 
 Finally, the term in \eqref{A22} is the one we are estimating and will be handled by Gronwall.
 
 Hence, gathering the above information we conclude in the case $j=8$ that
 \begin{equation}\label{A26}
 \begin{split}
 & \underset{[0,T]}{\sup} \int (\p_x^8u)^2(x,t)\,\psi_{8,v, \epsilon, b}(x,t)\,dx\\
 &\;\;\; +\int\limits_0^T\int (\p_x^9u)^2(x,t)\, \chi_{_{\epsilon, b}}'(x+vt)\,dxdt\\
&  \le c = c(\|u_0\|_{7,2}; \|\p_x^8u_0\|_{L^2((0,\infty))}; v; \epsilon; b).
  \end{split}
  \end{equation}

We notice that by \eqref{A15} that \eqref{A26} implies that  for $\epsilon>0$,
$b>5\epsilon$, $v>0$,
\begin{equation}\label{A27}
\begin{split} 
 & \underset{[0,T]}{\sup} \int (\p_x^8u)^2(x,t)\,\chi_{_{ \epsilon, b}}(x+vt)\,dx\\
 &\;\;\; +\int\limits_0^T\int (\p_x^9u)^2(x,t)\, \chi_{_{\epsilon, b}}'(x+vt)\,dxdt\\
&  \le c = c(\|u_0\|_{7,2}; \|\p_x^8u_0\|_{L^2((0,\infty))}; v; \epsilon; b).
  \end{split}
  \end{equation}

Once we have established the desired result \eqref{A27} for the case $j=8$ we
sketch the iterative argument for the general case $j=8, \dots, m$.

Assuming the step $j=m_0\in \{ 8, \dots, m\}$, i.e. for $j=8, \dots,m_0$,
\begin{equation}\label{A28}
\begin{split} 
 & \underset{[0,T]}{\sup} \int (\p_x^j u)^2(x,t)\,\chi_{_{ \epsilon, b}}(x+vt)\,dx\\
 &\;\;\; +\int\limits_0^T\int (\p_x^{j+1}u)^2(x,t)\, \chi_{_{\epsilon, b}}'(x+vt)\,dxdt\\
&  \le c = c(\|u_0\|_{7,2}; \|\p_x^{m_0}u_0\|_{L^2((0,\infty))}; v; \epsilon; b),
  \end{split}
  \end{equation}
we shall prove it for $j=m_0+1$.

We repeat the argument in \eqref{A7}-\eqref{A8} to get $E_{1,m_0+1}$,
$E_{2, m_0+1}$, and $E_{3, m_0+1}$.

The estimate for  $E_{2, m_0+1}$ is similar to that given in \eqref{A9}-\eqref{A14}.

To handle the term  $E_{1,m_0+1}$ we observe that from \eqref{A20} (with $j=m_0+1$) one
has an estimate as in \eqref{A20}-\eqref{A20b} with $m_0+1$ instead of $8$, i.e. the
terms $E_{1, m_0+1,1}$ and $E_{1,m_0+1,2}$. As in the previous case $E_{1, m_0+1,1}$
is a multiple of the term we are estimating and $E_{1, m_0+1,2}$  can be bounded, 
after time integration, by the second term in the right hand side of \eqref{A28} by taking
there an appropriate value of $\epsilon$ and $b$. So it only remains to consider
the terms in $E_{3, m_0+1}$, where (see \eqref{A8})
\begin{equation*}
E_{3, m_0+1}= \int Q_{m_0+1} \big( (\p_x^lu)^2_{|l|\le m_0+1} ; (\p_x^r\psi_{m_0+1})_{|r|\le3}\big)
\,dx.
\end{equation*}
where $Q_{m_0+1}(\cdot,\cdot) $ is a polynomial in its variables, linear in the 
($\psi_{m_0+1}$, $\p_x\psi_{m_0+1}$, $\p_x^2\psi_{m_0+1}$, $\p_x^3\psi_{m_0+1}$)
components and at most quadratic in the highest derivatives of $u$, i.e. $(\p_x^{m_0+1}u)$ involving at most $7+2(m_0+1)$ derivatives of $u$.

To handle $E_{3,m_0+1}$ one combines the global (in space) estimate in \eqref{A1}
with that in \eqref{A0b} and those obtained in the previous estimates $j=8,\dots,m_0$,
i.e. \eqref{A28} with $j=8,\dots, m_0$ to obtain the desired estimate. The proof follows
the argument provided in details in \cite{ILP-cpde} Section 3. Therefore it will be omitted.

Finally, to justify the previous formal computation we recall that the argument in the proof of  Theorem \ref{thm1} shows that $u$ in \eqref{A1}-\eqref{A0b}
is the limit in the $C([0,T] : H^{7-\delta}(\R))$-norm (for any $\delta>0$) of smooth 
solutions  (weak form of the continuous dependence upon the data). In particular, we
have that $u$ is the uniform limit of smooth solutions in $\R\times[0,T]$. Hence,
by performing the above (formal) argument in the smooth solutions one obtains
a uniform bounded sequence in the norms described in \eqref{i3} and \eqref{i4}. Hence,
considering the uniform boundedness, the  weak convergence and passing to the limit we obtain the desired result.

\section{Dispersive Blow-up}

Consider the IVP associated to the $k$-generalized Korteweg-de Vries equation,
\begin{equation}\label{gmkdv}
\begin{cases}
\partial_tu+\partial_x^3u +u^k\partial_xu=0, \hskip5pt k=1, 2, 3, \dots, \;x\in\R, \;t>0,\\
u(x,0)=u_0(x).
\end{cases}
\end{equation}

\vspace{3mm}
\noindent\underline{\bf Proof of Theorem \ref{thm3}.}
\vspace{3mm}

Let $\phi(x)=e^{-2|x|}$ and consider  the linear IVP,
\begin{equation}\label{linear-kdv}
\begin{cases}
\partial_t v+\partial_x^3v=0, \hskip 5pt x\in\R, \;t>0,\\
v(x,0)=v_0
\end{cases}
\end{equation}
whose solution is given by
\begin{equation}\label{1.4b}
v(x,t)= V(t)v_0(x)= e^{-t\partial_x^3}v_0= S_t\ast v_0(x)
\end{equation}
where
\begin{equation*}
S_t(x)=\frac{1}{\sqrt[3]{3t}} A_i\Big(\frac{x}{\sqrt[3]{3t}}\Big)
\end{equation*}
and $A_i(\cdot)$ denotes the Airy function.

Define
\begin{equation}\label{data}
u_0(x)=\underset{j=1}{\overset{\infty}{\sum}} \alpha_j\,V(-j)\phi(x), \hskip10pt \alpha_j>0.
\end{equation}

If 
\begin{equation*}
\underset{j=1}{\overset{\infty}{\sum}}  \alpha_j \ll 1,
\end{equation*}
we have that $u_0\in H^1(\R)$, in fact in $u_0\in H^{3/2-}(\R)$. This in particular guarantees the global existence of solutions in $H^1(\R)$ for the IVP \eqref{gmkdv}.\\

\noindent\underline{\bf Step 1.}{ Reduction to linear case.} 
\vskip3mm
\noindent \underline{\bf Case $k=2, 3, \dots$}
\vskip2mm

Since the nonlinear part of the solution $u$ of \eqref{gmkdv}, i.e.
\begin{equation*}
z(t)=\int\limits_0^t V(t-t') \,u^k\partial_xu(t')\,dt'
\end{equation*}
is in $C([0,\infty) : H^2(\R))$ (see \cite{LS} for the proof in the case  $k=2$, since our solution is  $C([0,T]: H^{1}(\R))$, the argument works for all $k\ge2$). It suffices to consider the linear part, $V(t)u_0$.\\

\noindent \underline{\bf Case $k=1$}
\vskip2mm

We observe that from \eqref{1.4b}--\eqref{data}  it follows that $u_0\in H^{{3/2}^{-}}(\R)$.

Next we recall the identity deduced in \cite{FLP}: for $\beta\in (0,1)$ and $t\in\R$ 
\begin{equation}\label{weight1}
|x|^{\beta} V(t)f = V(t)(|x|^{\beta}f)+ V(t)\{\Phi_{t,\beta}(\widehat{f\,})\}^{\vee}
\end{equation}
with
\begin{equation}\label{weight2}
\|\Phi_{t, \beta}(\widehat{f\,})\|_2\le c(1+|t|) \|f\|_{2\beta,2}.
\end{equation}
Hence, for $\beta\in (0,3/4)$
\begin{equation*}
|x|^{\beta}u_0\in L^2(\R) \hskip20pt \text{if} \hskip20pt \underset{j=1}{\overset{\infty}{\sum}} \,\alpha_j \,j^{\beta}<\infty.
\end{equation*}

Assuming the last inequality we have that for $u_0$ as above the corresponding solutions of the IVP for the KdV equation satisfies
for any $T>0$
\begin{equation}
u\in C([0,T] : H^{{3/2}^{-}}(\R)\cap L^2(|x|^{{3/4}-\epsilon}\,dx)),
\end{equation}
\begin{equation}
J^{3/2^{-}}\p_xu\in L^{\infty}_x(\R : L^2([0,T])), \hskip15pt \text{(smoothing effect)}
\end{equation}
and
\begin{equation}\label{strichartz-kdv}
\begin{split}
\int\limits_0^T \|&D^{\alpha\theta/2} J^s u(\cdot, t)\|_p^q\,dt<\infty  \hskip7pt\text{for}\;\;0<s<3/2, \hskip15pt \text{(Strichartz)}\\
&(q,p)=(6/\theta(\alpha+1), 2/(1-\theta)), \;\;\theta\in(0,1),\;\; 0\le \alpha\le 1/2.
\end{split}
\end{equation}

As in the previous case $k=2,3,\dots$ we shall show that for any $t\in [0,T]$
\begin{equation}
z(t)=\int\limits_0^t V(t-t') u\p_x u (t')\,dt' \in C^1(\R)
\end{equation}
by proving that
\begin{equation}
z(t)\in C([0,T] : H^{3/2^{-}+1/6}(\R)).
\end{equation}

Using the inequality (see \cite{KPV93})
\begin{equation}
\underset{0\le t\le T}{\sup}\|\p_x \int\limits_0^t V(t-t') F(\cdot, t')\,dt'\|_2 \le c\|F\|_{L^1_xL^2_T}
\end{equation}
one has that
\begin{equation*}
\begin{split}
\underset{0\le t\le T}{\sup}\|D^{3/2^{-}+1/6}&\int\limits_0^t V(t-t') u\p_xu (t')\, dt'\|_2 \le \|D^{1/2^{-}+1/6}(u\p_x u)\|_{L^1_xL^2_T}\\
&\le \big( \|u\|_{L^{6/5}_xL^3_T}\|D^{3/2^{-}+1/6}u\|_{L^6_xL^6_T} +E_1\big)
\end{split}
\end{equation*}
where the terms in $E_1$ are easy to control by considering the commutator estimates in the Appendix of \cite{KPV93}  and interpolated
norms of the previous terms to be considered below, so we omit this proof.

Now
\begin{equation*}
\|D^{3/2^{-}+1/6}u\|_{L^6_xL^6_T}<\infty
\end{equation*}
from \eqref{strichartz-kdv} with $p=q=6$, $\theta=2/3$, $\alpha=1/2$ and using the inequality (2.11) in \cite{FP},
\begin{equation*}
\|J^{\gamma\,a}(\langle x\rangle^{(1-\gamma) b}f)\|_2 \le c\|\langle x\rangle^b f\|_2^{1-\gamma}\|J^af\|_2^{\gamma},\;\;a,b>0, \;\gamma\in(0,1),
\end{equation*}
we deduce
\begin{equation*}
\begin{split}
 \|u\|_{L^{6/5}_xL^3_T}&\le c\|\langle x\rangle^{1/2+} u\|_{L^3_TL^3_x}\le c T^{1/3}\|\langle x\rangle^{1/2+} u\|_{L^{\infty}_TL^3_x}\\
 &\le c T^{1/3}\|J^{1/6} (\langle x\rangle^{1/2+} u)\|_{L^{\infty}_TL^2_x}\\
 &\le cT^{1/3}\|J^{3/2-}u\|_{L^{\infty}_TL^2_x}^{1-\gamma}\|\langle x\rangle^{3/4-} u\|_{L^{\infty}_TL^2_x}^{\gamma}
 \end{split}
 \end{equation*}
 with $\gamma$ such that $\gamma \,{3/4}^{-}= {1/2}^{+}$ (i.e. $\gamma>2/3$) and $(1-\gamma){3/2}^{-}>1/6$.

As in the case $k=2, 3, \dots$ we have reduced ourselves to consider the linear associated problem so the nonlinearity after Step 1 is not relevant for our purposes.

\vskip5mm
\noindent\underline{\bf Step 2.} Estimate for $V(t)\phi$, $(t>0)$. 

Assume that $v_0\in L^2(\R)$ and $e^{x}v_0\in L^2(\R)$

Now consider $w(x,t)=e^x v(x,t)$. Following Kato \cite{Kato} we set  $v(x,t)= e^{-x} w(x,t)$ where $w$ is solution of 
\begin{equation}\label{ivpw}
\begin{cases}
\partial_t w +(\partial_x-1)^3 w=0,\\
w(x,0)= e^x v_0(x).
\end{cases}
\end{equation}
Since
\begin{equation*}
(\partial_x - 1)^3= \partial_x^3 -3\partial_x^2+3\partial_x-1
\end{equation*}
one has that
\begin{equation*}
\begin{split}
w(x,t)&= V(t) e^{3t\partial_x^2} e^{-3t\partial_x} e^t \big( e^x v_0(x)\big)\\
&= V(t) e^{3t\partial_x^2} e^{-3t\partial_x} \big(e^{x+t} v_0(x)\big)
\end{split}
\end{equation*}
and
\begin{equation*}
V(t)v_0= v(x,t)= e^{-x} V(t)e^{3t\partial_x^2} \big(e^{x-2t} v_0(x-3t)\big).
\end{equation*}

We notice using the heat kernel properties that
\begin{equation*}
\partial_x^m V(t)v_0 \sim e^{-x} V(t)\Big(\partial_x^m e^{3t\partial_x^2}\Big)\Big( e^{x-2t} v_0(x-3t)\Big).
\end{equation*}

It follows that
\begin{equation*}
\|\partial_x^m V(t)v_0\|_2\sim \frac{c_m}{(3t)^{m/2}}\| e^{x-2t} v_0(x-3t)\|_2\sim  \frac{c_m}{(3t)^{m/2}} e^t,
\end{equation*}
since
\begin{equation*}
\|e^{x-2t} v_0(x-3t)\|_2= e^t\|e^{x-3t} v_0(x-3t)\|_2= c\,e^t.
\end{equation*}

Similarly, if $t<0$, we have an IVP analogous to the one in \eqref{ivpw}   the operator $-(\partial_x+1)^3$ instead of $(\partial_x-1)^3$.
Thus
\begin{equation*}
\begin{split}
V(t)v_0&= e^x V(t) e^{-3t\partial_x^2} e^{-3t\partial_x} (e^{-t} e^{-x} v_0)\\
&= e^x V(t) e^{-3t\partial_x^2} e^{-3t\partial_x} \Big(e^{-x-t} v_0(x)\Big)\\
&= e^{x} V(t)  e^{-3t\partial_x^2} \Big( e^{-x-4t} v_0(x-3t)\Big).
\end{split}
\end{equation*}
and so  we have
\begin{equation*}
\partial_x^m V(t)v_0 \sim e^{-x} V(t)\Big(\partial_x^m e^{-3t\partial_x^2}\Big)\Big( e^{x-4t} v_0(x+3t)\Big)
\end{equation*}
and
\begin{equation}\label{suma}
\|\partial_x^m V(t)v_0\|_2\sim  e^x\,\frac{c_m}{(3t)^{m/2}} e^{-t}.
\end{equation}
\\

\noindent\underline{\bf Step 3.} Next we prove that
\begin{equation*}
\underset{j=1}{\overset{\infty}{\sum}} \alpha_j V(-j) \phi \in C^{\infty}(\R)
\end{equation*}
or equivalently
\begin{equation*}
\underset{j=1}{\overset{\infty}{\sum}} \alpha_j e^x V(-j) \phi \in C^{\infty}(\R).
\end{equation*}

To do this, it suffices to show that 
\begin{equation*}
\underset{j=1}{\overset{\infty}{\sum}} \alpha_j e^x\Big(\partial_x^m V(-j) \phi\Big) \in 
L^2_{\rm loc}(\R) \;\; \text{ for all} \;\;m
\end{equation*}
or equivalently
\begin{equation*}
 \underset{j=1}{\overset{\infty}{\sum}} \alpha_j \frac{c_m}{(3j)^{m/2}}\, e^j <\infty.
 \end{equation*}
\\

\noindent\underline{\bf Step 4.} For each $t>0$, $t\notin \Z^{+}$, we claim that 

\begin{equation*}
V(t)u_0= \underset{j=1}{\overset{\infty}{\sum}} \alpha_j  V(t-j)\phi \in C^1(\R).
\end{equation*}

Combining  \eqref{suma}  and the assumption
\begin{equation*}
\underset{j=1}{\overset{\infty}{\sum}}\, \alpha_j  \frac{1}{3|t-j|} e^{|t-j|} <\infty
\end{equation*}
one has $V(t)u_0 \in H^2_{\rm loc} (\R)\subseteq C^1(\R)$.\\

\noindent\underline{\bf Step 5.}  For $t=n\in \Z^{+}$ we affirm that

\begin{equation}\label{suma2}
V(n)u_0 = \alpha_n\phi +\underset{j\neq n}{ \underset{j=1}{\overset{\infty}{\sum}} \alpha_j V(n-j) \phi}\equiv \alpha_n\phi +\Phi_n
\end{equation}
with $\Phi_n\in C^1$.

As before using \eqref{suma} and taking
\begin{equation*}
\underset{j\neq n}{ \underset{j=1}{\overset{\infty}{\sum}}}\, \alpha_j  \frac{1}{3|n-j|} e^{|n-j|} <\infty
\end{equation*}
it follows that $\Phi_n\in H^2_{\rm loc} (\R)$ which yields \eqref{suma2}.

By setting $\alpha_j= c\,e^{-j^2}$ with $c$ small enough we obtain the desired result.


\vspace{3mm}

\noindent\underline{\bf Proof of Theorem \ref{thm4}.}
\vspace{3mm}

 (a) First we consider the case $k=2, 3, \dots$. We recall the Strichartz estimates for solutions of the linear IVP \eqref{linear-kdv}
established in \cite{KPV91}
\begin{equation}\label{kdv-strichartz1}
\Big(\int\limits_{-\infty}^{\infty} \| D^{\alpha\theta/2} V(t)f\|_p^q\,dt\Big)^{1/q} \le c\|f\|_2,
\end{equation}
with $(q, p)= (6/{\theta(\alpha+1)}, 2/(1-\theta))$, $0\le \theta\le 1\;$ and $\;0\le \alpha\le 1/2$.

In particular for $p\in(2,\infty)$ and $\alpha=1/2$, the estimate
 \eqref{kdv-strichartz1}  becomes
\begin{equation}\label{kdv-strichartz2}
\Big(\int\limits_{-\infty}^{\infty} \| D^{(p-2)/4p} V(t)f\|_p^{4p/(p-2)}\,dt\Big)^{(p-2)/4p} \le c\|f\|_2.
\end{equation}

We take $\tilde{u}_0\in H^s(\R)$ with $s=j-\frac{p-2}{4p}=j-\hat{p}>3/4$ with 
$\tilde{u}_0\notin W^{j,p}(\R^{+})$.

From \eqref{kdv-strichartz2} it follows that there exists $\hat{t}\in (0,T/2)$ such that 
\begin{equation*}
V(\pm\hat{t}\,)\tilde{u}_0,\; V(\pm 2\hat{t}\,)\tilde{u}_0 \in W^{r,p}\hskip10pt\text{with}\hskip10pt r=s+\frac{p-2}{4p}=j.
\end{equation*}

Thus we consider the initial data
\begin{equation}
u_0=V(\hat{t}\,)\tilde{u}_0+ V(-\hat{t}\,)\tilde{u}_0. 
\end{equation}

Observe that  $u_0\in H^s(\R)$, so  since
\begin{equation*}
u(t)=V(t)u_0-\int\limits_0^t V(t-t') u^k \partial_x u (t')\,dt'= V(t)u_0 + z(t),
\end{equation*}
from the argument in \cite{LS} one has that
\begin{equation*}
z\in C([-T,T]:H^{s+1}(\R))  \hookrightarrow C([-T,T] : W^{j,p}(\R)).
\end{equation*}

Also one sees that
\begin{equation*}
V(\hat{t}\,)u_0= V(2\hat{t}\,)\tilde{u}_0+\tilde{u}_0 \notin W^{j,p}(\R^{+}).
\end{equation*}

Similarly for $V(-\hat{t}\,)u_0$. Gathering this information we obtain the desired result.\\

 (b) Now we turn to the proof of the case $k=1$. We observed  that the argument of proof in Theorem \ref{thm3} (Step 1) shows that in the case $k=1$
if $\tilde{u}_0\in H^{\hat{\!j}}(\R)\cap L^2(|x|^{j/2}\,dx)$ with $\hat{\!j}=j+1/2-1/p-1/{12}$ (thus $H^{\hat{\!j}}\hookrightarrow W^{j,p}$) with $\tilde{u}_0\notin W^{j,p}(\R)$
one has that
\begin{equation*}
z(t)=\int\limits_0^t V(t-t') u\p_xu(t')\,dt'\in C([-T,T] : H^{j+\frac12-\frac{1}{p}-\frac{1}{12}}(\R))\hookrightarrow W^{j,p}(\R).
\end{equation*}

Once this has been established the rest of the proof follows the argument provided for the case $k=2, 3, \dots$.

\section{BBM equation}
\noindent\underline{\bf Proof of Theorem \ref{thm5}.}
\vspace{3mm}

We shall restrict  ourselves  to consider the most general case $s=0$, i.e. $u_0\in L^2(\R)$.
Thus, from the local well-posedness theory  in Theorem A4 (\cite{BT})  there
exist $T=T(\|u_0\|_{2})>0$ and a unique solution $u=u(x,t)$ of the IVP \eqref{ivp-bbm} such that 
\begin{equation}
\label{A0}u\in C([0,T] : L^2(\R)).
\end{equation}
We rewrite the BBM equation in \eqref{ivp-bbm}
\begin{equation*}
\p_tu+\p_xu+u\p_xu-\p_x^2\p_tu=0
\end{equation*}
as the integro-differential equation  
\begin{equation}\label{bbm1}
\p_t u=-\p_x J^{-2}(u+u^2/2)
\end{equation}
where
\begin{equation*}
J^{-2}f=(1-\p_x^2)^{-1}f=\frac12 e^{-|x|}\ast f.
\end{equation*}
We observe that
\begin{equation}
\label{A1b}
\p_x^2J^{-2}=J^{-2}-I.
\end{equation}

Since  $u\in C([0,T] : L^2(\R))$ Sobolev embedding theorem guarantees that
\begin{equation}
\label{A2b}
\p_xJ^{-2} u \in C([0,T] : H^1(\R)) \hookrightarrow C([0,T] : C_{\infty}(\R)) ,
\end{equation}
where 
\begin{equation*}
\,C_{\infty}(\R)=\{ f :\R\to \R\,:\,f\;\text{continuous  with}\;\;\lim_{|x|\to\infty}f(x)=0\}.
\end{equation*}
Also, since
\begin{equation*}
u^2\in C([0,T] : L^{1}(\R)),
\end{equation*}
one has that 
\begin{equation}
\label{A3b}
\p_xJ^{-2} (u^2)=-  \frac12 \sgn(x) e^{-|x|}\ast u^2\in C([0,T]: C_b(\R)),
\end{equation}
where $\,C_{b}(\R)=C(\R)\cap L^{\infty}(\R)$. Hence, combining \eqref{bbm1}, \eqref{A2b} and \eqref{A3b} it follows that
\begin{equation}
\label{bbm2}
\begin{aligned}
u(x,t)&= u_0(x)-\int\limits_0^t \p_xJ^{-2}(u+u^2/2)(x,\tau)\,d\tau \\
&= u_0(x) + z(x,t),
\end{aligned}
\end{equation}
with
\begin{equation*}
z(x,t)=z\in C([0,T]: C_b(\R)).
\end{equation*}

Thus, if for some open $\Omega\subset \R$, $u_0\big|_{\Omega}\in C(\Omega)$, then
\begin{equation}\label{A5b}
u\big|_{\Omega\times[0,T]} \in C([0,T]:C(\Omega)).
\end{equation}
Moreover,  using \eqref{A1b}  one has that
\begin{equation}\label{A6b}
\begin{aligned}
\p_x z(x,t)&=-\int\limits_0^t \p^2_xJ^{-2}(u+u^2/2)(x,\tau)\,d\tau \\
&=\int\limits_0^t ((u+u^2/2)-J^{-2}(u+u^2/2))(x,\tau)\,d\tau. 
\end{aligned}
\end{equation}
Thus, from  \eqref{A5b}
\begin{equation*}
(u+u^2/2)\big|_{\Omega\times[0,T]} \in C([0,T]:C(\Omega)),
\end{equation*}
and by \eqref{A0} and an argument similar to that in \eqref{A2b}-\eqref{A3b} 
\begin{equation*}
J^{-2}(u+u^2/2)\in C([0,T]:C_b(\R)),
\end{equation*}
therefore we conclude that if  $u_0\big|_{\Omega}\in C(\Omega)$, then
\begin{equation*}
z\big|_{\Omega\times[0,T]}\in C([0,T]: C^1(\Omega)).
\end{equation*}

Hence, from \eqref{bbm2} :  if for some $\Omega \subset \R$  open, $u_0\big|_{\Omega}\in C^{\theta}(\Omega)$ with $\theta\in(0,1]$, then
\begin{equation}\label{A7a}
u\big|_{\Omega\times[0,T]} \in C([0,T]:C^{\theta}(\Omega)).
\end{equation}

Now using the previous step  \eqref{A7a} with $\theta=1$,  i.e. if for $\Omega \subset \R$ open, $u_0\big|_{\Omega}\in C^{1}(\Omega)$, then
 \begin{equation}\label{A7b}
u\big|_{\Omega\times[0,T]} \in C([0,T]:C^1(\Omega)).
\end{equation}
This combined with \eqref{A6b} implies that
\begin{equation*}
z(x,t)=z\in C([0,T]: C^2(\Omega)).
\end{equation*}
Therefore, if $u_0\big|_{\Omega}\in C^{1+\theta}(\Omega)$ for some $\theta\in (0,1]$, then from \eqref{bbm2} it follows that
 \begin{equation}\label{A8b}
u\big|_{\Omega\times[0,T]} \in C^1([0,T]:C^{1+\theta}(\Omega)).
\end{equation}

It is clear that by reapplying this argument one gets the desired result.

\section{Degasperis-Procesi equation}
In this section we shall consider the IVP associated to the Degasperis-Procesi (DP) equation
\begin{equation}\label{ivpDP}
\begin{cases}
\partial_tu-\partial_t\partial_x^2u+4u\partial_xu= 3\partial_xu\partial_x^2u+u\partial_x^3u,\hskip5pt x\in\R, \;t>0,\\
u(x,0)=u_0(x).
\end{cases}
\end{equation}

 The equation in \eqref{ivpDP} can be rewritten in the integro-differential form
\begin{equation}\label{DP}
\partial_tu+u\partial_x u+\frac32 \,(1-\partial_x^2)^{-1}\partial_x (u^2)= 0,
\end{equation}
where
\begin{equation*}
(1-\partial_x^2)^{-2}f=J^{-2}f= \frac12 e^{-|x|}\ast f.
\end{equation*}
 Notice that
 \begin{equation*}
 \partial_x^2 J^{-2}  = J^{-2} - I.
 \end{equation*}
 
 \vspace{3mm}
 
 \noindent\underline{\bf Proof of Theorem \ref{thm6}.}
\vspace{3mm}
 
In \cite{yin} it was shown that the IVP associated to the equation \eqref{DP} is locally well-posed in $H^s(\R)$ for $s>3/2$.

Let $u^{\epsilon}$ be the  solution corresponding to the initial data $\rho_{\epsilon}\ast u_0 =u_0^{\epsilon}$, with $u_0\in H^{1+\delta}(\R)\cap W^{1,\infty}(\R),\;\delta>0,$ and 
$\rho_{\epsilon}$ denoting the usual mollifiers. Thus, 
\begin{equation}
\label{B1}
u^{\epsilon}\in C([0,T_{\epsilon}]:H^{\infty}(\R))\cap...
\end{equation}
To estimate $T_{\epsilon}$ we recall that using the commutator estimates in \cite{KP}, see \eqref{kapo1}-\eqref{kapo2}, and simpler inequalities as those below one obtains the  formal energy estimate: for any $v_0\in H^s(\R)$ with $s> 3/2$ 
the corresponding solution to the DP equation  in \eqref{DP} $v\in C([0,T] : H^s(\R))$ with $T=T(\|v_0\|_{s,2})>0$ obtained in \cite{yin} satisfies :\begin{equation}
\label{B2}
\frac{d}{dt}\| v (t)\|_{s,2}\leq c_s (\|v(t)\|_{\infty} +\|\partial_x v(t)\|_{\infty})\,\| v(t)\|_{s,2}.
\end{equation}

The estimate \eqref{B2} implies  the following continuation principle : given $v_0\in H^s(\R),\,s>3/2$, then the  corresponding solution $v\in C([0,T] : H^s(\R))$ of the IVP associated to \eqref{DP} can be extended 
in the time interval $[0,T^*]$ with $T^*>T$ satisfying that $v\in C([0,T^*]:H^s(\R))$ whenever 
\begin{equation}
\label{B3a}
\int_0^{T^*}(\|v(t)\|_{\infty} +\|\partial_x v(t)\|_{\infty})\,dt <\infty.
\end{equation}

\vskip5mm
\noindent\underline{\bf A priori estimate.}
\vskip2mm

We shall show that if $u_0\in H^{1+\delta}(\R)\cap W^{1,\infty}(\R),\,\delta>0$, then $T_{\epsilon} $ defined above can be estimated independently  of $\epsilon$, i.e. $T_{\epsilon}=T,\;\forall \epsilon>0$, with
\begin{equation*}
T=O( (\|u_0\|_{1,2}+\|u_0\|_{1,\infty})^{-1})\hskip12pt\text{as}\hskip12pt(\|u_0\|_{1,2}+\|u_0\|_{1,\infty})\downarrow 0.
\end{equation*} 

Applying energy estimates we have that
\begin{equation}
\frac{d}{dt} \|u^{\epsilon}(t)\|_2\le c\,\|u^{\epsilon}(t)\|_{\infty}\|u^{\epsilon}(t)\|_2
\end{equation}
after using that 
\begin{equation*}
\|\partial_x J^{-2}(u^{\epsilon})^2\|_2=c \| \,\sgn(x)e^{-|x|}\ast(u^{\epsilon})^2\|_2 \le \|(u^{\epsilon})^2\|_2\le \|u^{\epsilon}\|_{\infty}\|u^{\epsilon}\|_2.
\end{equation*}

Also, as far as the characteristics flow is  defined, i.e. $X^{\epsilon}(t; x_0)$,
\begin{equation*}
\begin{cases}
\dfrac{dX^{\epsilon}}{dt}(t) = u^{\epsilon} ( X^{\epsilon}(t), t),\\
\\
X^{\epsilon}(0)=x_0,
\end{cases}
\end{equation*}
one has that
\begin{equation*}
\frac{d}{dt} u^{\epsilon}(X^{\epsilon}(t; x_0), t)= -\frac32 \partial_x J^{-2} ( (u^{\epsilon})^2).
\end{equation*}

Since
\begin{equation*}
\|\partial_x J^{-2}(u^{\epsilon})^2\|_{\infty} \le \|(u^{\epsilon})^2\|_1\le \|u^{\epsilon}\|_2^2,
\end{equation*}
it follows that
\begin{equation}
\frac{d}{dt} \|u^{\epsilon}(t)\|_{\infty} \le c\,\|u^{\epsilon}(t)\|_2^2.
\end{equation}

Next, we recall the following estimates deduced in \cite{KP}: for any $r>0$
\begin{equation}
\label{kapo1}
\| [J^r;f]g\|_2\leq c_r(\|\partial_xf\|_{\infty}\|J^{r-1}g\|_2+\|g\|_{\infty}\|J^rf\|_2),
\end{equation}
and
\begin{equation}
\label{kapo2}
\| J^r(fg)\|_2\leq c_r(\|f\|_{\infty}\|J^{r}g\|_2+\|g\|_{\infty}\|J^rf\|_2).
\end{equation}
Combining \eqref{kapo1}-\eqref{kapo2} one gets that
\begin{equation}
\begin{split}
\frac{d}{dt}\|u^{\epsilon}&(t)\|_{1+\delta,2}\\
&\le c(\|\partial_x u^{\epsilon}(t)\|_{\infty}\|u^{\epsilon}(t)\|_{\delta,2}+ c\|u^{\epsilon}(t)\|_{\infty}\|u^{\epsilon}(t)\|_{1+\delta,2}).
\end{split}
\end{equation}
Finally, since
\begin{equation}
\partial_t(\partial_x u^{\epsilon})+u^{\epsilon}\partial_x(\partial_x u^{\epsilon})+(\partial_x u^{\epsilon})(\partial_x u^{\epsilon})+
\frac32\Big(J^{-1} (u^{\epsilon})^2-(u^{\epsilon})^2\Big)=0,
\end{equation}
then
\begin{equation*}
\begin{aligned}
\frac{d}{dt}&(\partial_xu^{\epsilon}(X^{\epsilon}(t;x_0),t)) \\
&+\big((\partial_xu^{\epsilon})^2 +\frac{3}{2} (J^{-1}(u^{\epsilon})^2-(u^{\epsilon})^2)\big)(X^{\epsilon}(t;x_0),t)=0.
\end{aligned}
\end{equation*}
Thus, using that
\begin{equation*}
\|J^{-2} (u^{\epsilon})^2-(u^{\epsilon})^2\|_{\infty}\le c\,\|u^{\epsilon}\|_{\infty}^2,
\end{equation*}
we conclude that if $u_0\in H^{1+\delta}(\R)\cap W^{1,\infty}(\R)$ (actually $H^s(\R)$ with $s>1/2$ instead of $1+\delta$ will suffice for this step), there exists
$T=T(\|u_0\|_{1+\delta,2}; \|u_0\|_{1,\infty})>0$ such that
\begin{equation}\label{DP1}
\underset{[0,T]}{\sup} \Big(\|u^{\epsilon}(t)\|_{1+\delta,2}+\|u^{\epsilon}(t)\|_{1,\infty}\Big)
\le M\Big(\|u_0\|_{1+\delta,2}+\|u_0\|_{1,\infty}\Big).
\end{equation}
We recall that if $v\in C(\R\times [0,T])\cap L^{\infty}([0,T]:W^{1,\infty}(\R))$, then the characteristic flow $X_v(t; x_0)=X(t; x_0)$ solution of
\begin{equation}\label{flow-dp}
\begin{cases}
\dfrac{dX}{dt}(t) = v ( X(t), t),\\
\\
X(0)=x_0,
\end{cases}
\end{equation}
is well-defined.
Thus, combining these facts with the  continuation principle in \eqref{B2}-\eqref{B3a} we conclude that
\begin{equation*}
(u^{\epsilon})_{\epsilon>0}\subset C([0,T]:H^{\infty}(\R)),
\end{equation*}
with $T$ as in \eqref{DP1}. As it was remarked above for this step one just needs $u_0\in H^s(\R)\cap W^{1,\infty}(\R)$, $s>1/2$.
\vskip5mm

\noindent\underline{\bf Convergence as $\epsilon\downarrow 0$}
\vskip2mm

Defining $w=u^{\epsilon}-u^{\epsilon'}$ one gets the equation
\begin{equation}
\partial_t w +u^{\epsilon}\partial_x w+w\partial_x u^{\epsilon}+\frac32 \partial_x J^{-2}\big((u^{\epsilon}+u^{\epsilon'})w\big)=0.
\end{equation}
Thus
\begin{equation}\label{apriori-1}
\begin{split}											
\frac{d}{dt}\|w(t)\|_2 &\le c\Big( \|\partial_x u^{\epsilon}\|_{\infty}+\|\partial_x u^{\epsilon'}\|_{\infty}\Big)\|w(t)\|_2\\
&\;\;\; + c\Big( \|u^{\epsilon}\|_{\infty}+\| u^{\epsilon'}\|_{\infty}\Big)\|w(t)\|_2.
\end{split}
\end{equation}
Hence
\begin{equation*}
u^{\epsilon}\to u \hskip7pt \text{in} \hskip7pt C([0,T] :L^2(\R))\;\;\;\;\;\;\;\;\text{as}\;\;\;\;\;\;\;\epsilon\downarrow 0,
\end{equation*}
and consequently from \eqref{DP1}
\begin{equation}
\label{AA1}
u^{\epsilon}\to u \hskip7pt \text{in} \hskip7pt C([0,T] : H^{1}(\R)).
\end{equation}
Moreover, 
\begin{equation}
\label{AA2}
u\in C([0,T] : H^{1}(\R))\cap L^{\infty}([0,T]: H^{1+\delta}(\R)), 
\end{equation}
with
\begin{equation*}
\partial_t u^{\epsilon} \to \partial_t u\hskip7pt \text{in} \hskip7pt C([0,T] :L^2(\R))\;\;\;\;\;\;\;\;\text{as}\;\;\;\;\;\;\;\epsilon\downarrow 0.
\end{equation*}
This tells us that  $u=u(x,t)$ is the solution of the DP equation \eqref{DP} with data $u(x,0)=u_0(x)\in H^{1+\delta}(\R)\cap W^{1,\infty}(\R)$, where the equation is realized in $C([0,T]:L^2(\R))$.
Furthermore, since $u\in C([0,T] : C_b(\R))$  \eqref{DP1} and \eqref{AA1}  imply that $u\in L^{\infty}([0,T] :  W^{1,\infty} (\R))$ with norm  bounded by $M$   as in \eqref{DP1},
thus the characteristic flow $X_u(t;x_0)$, see \eqref{flow-dp}, is defined.

Notice that \eqref{apriori-1} implies uniqueness and a weak continuous dependence of the solutions upon the data, i.e. if $u_{0_n} \to u_0$ in $H^1$ with $u_{0_n}, u_0\in H^{1+\delta}(\R)$ uniformly bounded, then the corresponding solutions $u_n$ converges to $u$ in the $C([0,T]:H^1(\R))$-norm.

It is clear form our proof above that a weaker version of the Theorem \ref{thm6} holds for  $u_0\in H^{1+\delta}(\R)\cap W^{1,\infty}(\R)$, with $s>1/2$. However
we fixed $s=1+\delta$, $\delta>0$, such that the equation is realized in $C([0,T]: L^2(\R))$.

\vskip5mm

\noindent\underline{\bf Proof of Theorem \ref{thm7}.} 
\vspace{3mm}

Here we establish the propagation of regularity at the $C^1$ level
since at the $C^2$ level and beyond  it follows by writing $m=(1-\partial_x^2)u$ and considering the equation
\begin{equation*}
\partial_t m+u\partial_xm+3\partial_x u\, m=0.
\end{equation*}

Notice first that if $f=\frac32 \partial_x J^{-2}(u^2)$ with $u$ as in \eqref{AA2} one has that 
\begin{equation}
\label{AA3}
f \in C([0,T] : H^{2}(\R)). 
\end{equation}
Therefore, if 
$u_0\Big|_{\Omega}\in C^1$  with $\Omega \subseteq \R$ open,
then 
\begin{equation*}
u(\cdot,t)\Big|_{\Omega_t}\in C^1, \hskip10pt \Omega_t=\Phi_t(\Omega)
\end{equation*}
where $\Phi_t(x_0)=X_u(t; x_0)=X(t;x_0)$ defined as the solution of
\begin{equation*}
\begin{cases}
\dfrac{dX}{dt}=u(X(t), t),\\
\\
X(0)=x_0.
\end{cases}
\end{equation*}
Since the solution $u(\cdot, \cdot)$ satisfies
\begin{equation}
\begin{cases}
\partial_t u+u\partial_x u+f=0\\
u(x,0)=u_0,
\end{cases}
\end{equation}
with $f=f(x,t)$ as in \eqref{AA3} one sees that
\begin{equation*}
\frac{d}{dt}  u(X(t;x_0), t) = f(X(t; x_0), t)
\end{equation*}
and
\begin{equation*}
u(X(t;x_0), t)= u_0(x) +\int\limits_0^t f(X(\tau; x_0), \tau)\,d\tau.
\end{equation*}

This yields the desired result.


\section{Brinkman model 1--D case}

This section is concerned with the IVP associated to the Brinkman model,
\begin{equation}\label{brink}
\begin{cases}
\partial_t \rho=\partial_x \big(\rho(1-\partial_x^2)^{-1}\partial_x (\rho^2)\big), \;\;\;x\in \R, \;t>0,\\
\rho(x,0)=\rho_0(x).
\end{cases}
\end{equation}

 We shall use that 
\begin{equation}\label{brink2}
J^{-2}f= (1-\partial_x^2)^{-1}f =\frac12 e^{-|x|}\ast f
\end{equation}
and
\begin{equation}\label{brink3}
\partial_x^2 J^{-2} =J^{-2}  - I.
\end{equation}

\vskip.1in

\noindent\underline{\bf Proof of Theorem \ref{thm8}.}

\vspace{3mm}

Let $\rho^{\epsilon}$ be the solution corresponding to initial data $\rho_0^{\epsilon}(x)=G_{\epsilon}\ast \rho_0(x)$, with $G_{\epsilon}(x)=\epsilon^{-1}G\big(x/\epsilon\big)$, $G\in C^{\infty}_0(\R)$, $ G(x)\geq 0$, $\int G(x)dx =1$, and $\int x\,G(x)dx=0$. We recall that  in \cite{a-iorio}  the local well-posedness of the IVP 
\eqref{brink}  in $H^s(\R)$,
$s>3/2$ was established.

\vspace{3mm}
\noindent\underline{\bf A priori estimate in $H^1(\R)$.} 
\vspace{2mm}

To simplify the notation we shall use $\rho$ instead of $\rho^{\epsilon}$ in \eqref{B1}, \eqref{B2} and  \eqref{B5}. 

Energy estimates show that
\begin{equation}\label{B1b}
\begin{split}
\frac{d}{dt}\|\rho(t)\|_2^2&= \int \partial_x\rho (1-\partial_x^2)^{-1}\partial_x(\rho^2)\rho\,dx\\
&\;\;\;+
\int \rho (1-\partial_x^2)^{-1}\partial_x^2(\rho^2)\rho\,dx\\
&=\frac{1}{2}\,\int \rho (1-\partial_x^2)^{-1}\partial_x^2(\rho^2)\rho\,dx\\
& \le\frac{1}{2}  \|(J^{-2}  - I)(\rho^2)\|_{\infty}\|\rho(t)\|_2^2
\\
& \le c \|\rho\|_{\infty}^2\|\rho(t)\|_2^2.
\end{split}
\end{equation}
since $\|J^{-2}(\rho^2)\|_{\infty}\le c \|\rho^2\|_{\infty}$, and

\begin{equation}\label{B2b}
\begin{split}
\frac{d}{dt}\|\partial_x\rho(t)\|_2^2 
=& \int \partial_x^2\big(\rho(1-\partial_x^2)^{-1}\partial_x(\rho^2))\partial_x\rho\,dx\\
=& \int \partial_x^2\rho(1-\partial_x^2)^{-1}\partial_x(\rho^2)\partial_x\rho\,dx\\
&+2\int  \partial_x \rho(1-\partial_x^2)^{-1}\partial_x^2(\rho^2)\partial_x\rho\,dx\\
&+\int  \rho(1-\partial_x^2)^{-1}\partial_x^3(\rho^2)\partial_x\rho\,dx\\
\le &\, c\|\rho\|_{\infty}^2\|\partial_x\rho(t)\|_2^2.
\end{split}
\end{equation}

Notice that from Sobolev embedding theorem and \eqref{B1b}-\eqref{B2b} one gets that
\begin{equation*}
\frac{d}{dt}\|\rho^{\epsilon}(t)\|_{1,2} \leq c \|\rho^{\epsilon}(t)\|_{1,2}^3.
\end{equation*}
Therefore, there exists 
\begin{equation}\label{B3}
T=T(\|\rho_0^{\epsilon}\|_{1,2})=T(\|\rho_0\|_{1,2})=c \,\|\rho_0\|_{1,2}^{-2},
\end{equation}
 such that
\begin{equation}
\label{B4}
\sup_{[0,T]}\|\rho^{\epsilon}(t)\|_{1,2} \leq 2 \|\rho_0\|_{1,2}.
\end{equation}

Similarly,
\begin{equation}
\label{B5}
\frac{d}{dt}\|\partial_x^2\rho(t)\|_2 \le c\|\rho\|^2_{\infty}\|\partial_x^2\rho(t)\|_2.
\end{equation}

\vspace{3mm}
Combining \eqref{B4}-\eqref{B5} it follows that if $\rho_0\in H^1(\R)$, then
\begin{equation}
\label{B6}
\underset{[0,T]}{\sup} \|\partial^2_x\rho^{\epsilon}(t)\|_2= {\rm O}(\epsilon^{-1}).
\end{equation}

\vspace{10mm}
\noindent\underline{\bf Convergence as $\epsilon\downarrow 0$.}  
\vspace{2mm}

Let $\rho$ and $\tilde{\rho}$ be solutions of \eqref{brink}.  Thus, $w=\rho-\tilde{\rho}$ satisfies the equation
\begin{equation}\label{brink4}
\partial_t w-\partial_x(w(1-\partial_x^2)^{-1}\partial_x(\rho^2))-\partial_x(\tilde\rho (1-\partial_x^2)^{-1}\partial_x(\rho^2-{\tilde\rho}^2))=0.
\end{equation}

Multiplying \eqref{brink4} by $w$ and integrating in $x$ we have
\begin{equation}\label{brink5}
\begin{split}
\frac{d}{dt}\|w(t)\|_2^2&= \int \partial_x(w(1-\partial_x^2)^{-1}\partial_x(\rho^2))w\,dx\\
&\;\;\; +\int \partial_x(\tilde\rho (1-\partial_x^2)^{-1}\partial_x(\rho^2-{\tilde\rho}^2))w\,dx\\
&= W_1 + W_2.
\end{split}
\end{equation}
By integration by parts
\begin{equation}\label{brink6}
\begin{split}
W_1 &= \int  \partial_xw(1-\partial_x^2)^{-1}\partial_x(\rho^2) w\,dx\\
&\;\;\;\;+
\int w(1-\partial_x^2)^{-1}\partial_x^2(\rho^2)w\,dx\\
&\le c\,\|\rho^2\|_{\infty}\|w(t)\|_2^2
\end{split}
\end{equation}
and
\begin{equation}\label{brink7}
\begin{split}
W_2 &=\int \partial_x{\tilde\rho} (1-\partial_x^2)^{-1}\partial_x((\rho+{\tilde\rho})w)w\,dx\\
&\;\;\;+\int \tilde\rho (1-\partial_x^2)^{-1}\partial_x^2((\rho+{\tilde\rho}) w))w\,dx \\
& \le \|\partial_x{\tilde\rho}\|_2\|(1-\partial_x^2)^{-1}\partial_x((\rho+{\tilde\rho})w)\|_{\infty}\|w\|_2\\
&\;\;\;+\|\tilde{\rho}\|_{\infty}\|(\rho+\tilde{\rho})w\|_2\|w\|_2.
\end{split}
\end{equation}

Since
\begin{equation}\label{brink8}
\begin{split}
\|(1-\partial_x^2)^{-1}\partial_x((\rho+{\tilde\rho})w)\|_{\infty}\le \|(\rho+{\tilde\rho})w)\|_{1}\le 
\|\rho+{\tilde\rho}\|_2\|w\|_2,
\end{split}
\end{equation}
combining \eqref{brink5}, \eqref{brink6}, \eqref{brink7} and \eqref{brink8} we obtain  that
\begin{equation*}
\frac{d}{dt}\|w(t)\|_2^2 \le c (\underset{[0,T]}{\sup} \|\rho(t)\|_{1,2};\underset{[0,T]}{\sup}  \|\tilde{\rho}(t)\|_{1,2})\|w(t)\|_2^2.
\end{equation*}

Thus, if $0<\epsilon'<\epsilon$, then
\begin{equation}
\label{B7}
\underset{[0,T]}{\sup} \|(\rho^{\epsilon}-\rho^{\epsilon'})(t)\|_2 \le c(\|\rho_0\|_{1,2})
\|\rho^{\epsilon}_0-\rho^{\epsilon'}_0\|_2 ={\rm o}(\epsilon),\;\;\text{as}\;\;\epsilon\;\downarrow 0.
\end{equation}

Similarly, for any two strong solutions $\rho$, $\tilde{\rho}\in C([0,T]:H^1(\R))$  one has
\begin{equation*}
\underset{[0,T]}{\sup} \|(\rho-\tilde{\rho})(t)\|_2 \le c(\|\rho_0\|_{1,2};\|\tilde{\rho}_0\|_{1,2}) \,\|\rho_0-\tilde{\rho}_0\|_2.
\end{equation*}

Next, we shall estimate for $\|\partial_x(\rho^{\epsilon}-\rho^{\epsilon'})(t)\|_2$.  Let $w=\rho^{\epsilon}-\rho^{\epsilon'},$ so $w$ satisfies the equation
\begin{equation*}
\begin{split}
\partial_t w&-\partial_x(w(1-\partial_x^2)^{-1}\partial_x((\rho^{\epsilon})^2))\\
&-\partial_x(\rho^{\epsilon'} (1-\partial_x^2)^{-1}\partial_x((\rho^{\epsilon}+\rho^{\epsilon'})w))=0.
\end{split}
\end{equation*}
Thus
\begin{equation}
\label{B8}
\begin{split}
\frac{d}{dt} \|\partial_xw(t)\|_2^2&=\int \partial_x^2(w(1-\partial_x^2)^{-1}\partial_x((\rho^{\epsilon})^2))\partial_xw\,dx\\
&\;\;\;+\int \partial_x^2(\rho^{\epsilon'} (1-\partial_x^2)^{-1}\partial_x((\rho^{\epsilon}+\rho^{\epsilon'})w))\partial_x w\,dx\\
&= \widetilde{W}_1+\widetilde{W}_2.
\end{split}
\end{equation}

But
\begin{equation*}
\begin{split}
\widetilde{W}_1 &=\int \partial_x^2w(1-\partial_x^2)^{-1}\partial_x((\rho^{\epsilon})^2))\partial_xw\,dx\\
&\;\;\;+2\int \partial_x w(1-\partial_x^2)^{-1}\partial_x^2((\rho^{\epsilon})^2))\partial_xw\,dx\\
&\;\;\;+\int w(1-\partial_x^2)^{-1}\partial_x^3((\rho^{\epsilon})^2))\partial_xw\,dx
\end{split}
\end{equation*}
and so
\begin{equation*}
\begin{split}
|\widetilde{W}_1|\le c\|\rho^{\epsilon}\|_{\infty}^2 \|\partial w(t)\|_2^2+ c\|\rho^{\epsilon}\|_{1,2}^2
\|w(t)\|_2^{1/2}\|\partial_x w(t)\|_2^{3/2},
\end{split}
\end{equation*}
after integrating by parts and using that 
\begin{equation*}
\int w(1-\partial_x^2)^{-1}\partial_x^3((\rho^{\epsilon})^2)\partial_xw\,dx= - 
\int w(\partial_x-\partial_xJ^{-2})((\rho^{\epsilon})^2)\partial_xw\,dx.
\end{equation*}

Also,
\begin{equation*}
\begin{split}
\widetilde{W}_2 &= \int \partial_x^2\rho^{\epsilon'} (1-\partial_x^2)^{-1}\partial_x((\rho^{\epsilon}+\rho^{\epsilon'})w)\partial_x w\,dx\\
&\;\;\;+2\int \partial_x\rho^{\epsilon'} (1-\partial_x^2)^{-1}\partial_x^2((\rho^{\epsilon}+\rho^{\epsilon'})w)\partial_x w\,dx\\
&\;\;\;+\int \rho^{\epsilon'} (1-\partial_x^2)^{-1}\partial_x^3((\rho^{\epsilon}+\rho^{\epsilon'})w)\partial_x w\,dx\\
&= \widetilde{W}_{2,1}+\widetilde{W}_{2,2}+\widetilde{W}_{2,3}.
\end{split}
\end{equation*}

Then
\begin{equation*}
\begin{split}
 |\widetilde{W}_{2,1}|&\le \|\partial_x^2\rho^{\epsilon'}\|_2 \|(\rho^{\epsilon}+\rho^{\epsilon'})w\|_1\|\partial_x w\|_2\\
 &\le  \|\partial_x^2\rho^{\epsilon'}\|_2 \|\rho^{\epsilon}+\rho^{\epsilon'}\|_2\|w\|_2\|\partial_x w\|_2,
 \end{split}
 \end{equation*}
\begin{equation*}
\begin{split} 
|\widetilde{W}_{2,2}|&\le \|\partial_x\rho^{\epsilon'}\|_2 \|(\rho^{\epsilon}+\rho^{\epsilon'})w\|_{\infty}\|\partial_x w\|\\
&\le \|\partial_x\rho^{\epsilon'}\|_2\|\rho^{\epsilon}+\rho^{\epsilon'}\|_{\infty}\|w\|_{\infty}\|\partial_x w\|,
\end{split}
\end{equation*}
and  
\begin{equation*}
\begin{split}
|\widetilde{W}_{2,3}|&\le\|\rho^{\epsilon'}\|_{\infty} \|\partial_x((\rho^{\epsilon}+\rho^{\epsilon'})w)\|_2\|\partial_x w\|_2\\
&\le \|\rho^{\epsilon'}\|_{\infty}(\|\rho^{\epsilon}+\rho^{\epsilon'}\|_{\infty}\|\partial_xw\|_2+
\|\partial_x(\rho^{\epsilon}+\rho^{\epsilon'})\|_2\|w\|_{\infty})\|\partial_x w\|_2.
\end{split}
\end{equation*}

We observe that considering \eqref{B4}, \eqref{B6}, and \eqref{B7} one sees that
\begin{equation*}
 |\widetilde{W}_{2,1}| \le  c\, {\rm o}(1)\,\|\p_xw\|_2 \hskip10pt\text{as}\hskip10pt \epsilon \downarrow 0.
 \end{equation*}

Inserting the above estimates in \eqref{B8}, and then using \eqref{B6}-\eqref{B7} it follows that $\rho^{\epsilon}\to \rho$  with 
$\rho \in C([0,T] : H^1(\R))\cap C^1((0,T):L^2(\R))$ and  $T$ as in \eqref{B3}.
The continuous dependence of  the solution upon the data can be proved in a similar manner (see \cite{BoSm} or \cite{LP} Chapter 9). 
This basically completes the proof of Theorem \ref{thm8}

\vskip10pt
\noindent\underline{\bf  Proof of Theorem \ref{thm9}.}
\vskip5pt

Using \eqref{brink3} we write the equation in \eqref{brink} as
\begin{equation*}
\partial_t\rho-(1-\partial_x^2)^{-1}\partial_x(\rho^2)\partial_x\rho=(J^{-2}-I)(\rho^2)\rho.
\end{equation*}

Hence, formally one has that
\begin{equation}\label{br1}
\begin{aligned}
\partial_t\partial_x&\rho-(1-\partial_x^2)^{-1}\partial_x(\rho^2)\partial_x\partial_x \rho\\
&=2(J^{-2}-I)(\rho^2)\partial_x\rho-2\rho^2\partial_x\rho+\partial_xJ^{-2}(\rho^2)\rho,
\end{aligned}
\end{equation}

\begin{equation}\label{br2}
\begin{aligned}
\partial_t\partial^2_x&\rho-(1-\partial_x^2)^{-1}\partial_x(\rho^2)\partial_x\partial^2_x \rho\\
&=3(J^{-2}-I)(\rho^2)\partial^2_x\rho-2\rho^2\partial_x^2\rho+3\partial_x(J^{-2}-I)(\rho^2)\partial_x\rho\\
&\;\;\;+
\partial^2_xJ^{-2}(\rho^2)\rho-\big[\partial_x^2(\rho^2)\rho-2\rho^2\partial^2_x\rho],
\end{aligned}
\end{equation}

and for the general case $k\in\Z^+$

\begin{equation}\label{br3}
\begin{aligned}
\partial_t\partial^k_x&\rho-(1-\partial_x^2)^{-1}\partial_x(\rho^2)\partial_x\partial^k_x \rho\\
&=(k+1)(J^{-2}-I)(\rho^2)\partial^k_x\rho-2\rho^2\partial_x^k\rho\\
&\;\;\;+a^k_{k-1}\partial_x(J^{-2}-I)(\rho^2)\partial_x^{k-1}\rho\\
&\;\;\;+a^k_{k-2}\partial^2_x(J^{-2}-I)(\rho^2)\partial_x^{k-2}\rho\\
&\;\;\;+\dots+ a^k_{1}\partial^{k-1}_x(J^{-2}-I)(\rho^2)\partial_x\rho\\
&\;\;\;+\partial^k_xJ^{-2}(\rho^2)\rho-\big[\partial_x^k(\rho^2)\rho-2\rho^2\partial^k_x\rho].
\end{aligned}
\end{equation}

From Theorem \ref{thm8}  given $\rho_0\in H^1(\R)$, there exist $T>0$ (as in \eqref{B3}) and a unique strong solution $\rho\in C([0,T]: H^1(\R))$ of the IVP \eqref{brink}. We introduce the notation
\begin{equation}
\label{a1}
a(x,t)\equiv J^{-2}\partial_x(\rho^2)\in  C([0,T]: H^1(\R)) \hookrightarrow
C([0,T] :C^1_b(\R)),
\end{equation}
and
\begin{equation}
\label{a2}
(J^{-2}-I)\rho^2,\;\rho^2 \in C([0,T]: H^1(\R))\hookrightarrow
C([0,T] :C^1_b(\R)).
\end{equation}

In particular, from \eqref{a1} it follows that the flow $\Phi_t(x_0)=X(t;x_0)$ given by the solution of
\begin{equation}
\label{flow}
\begin{cases}
\begin{aligned}
&\frac{d\;}{dt}X=a(X,t),\\
\\
&X(0)=x_0,
\end{aligned}
\end{cases}
\end{equation}
is well defined for $t\in[0,T]$. For $\Omega\subset\R$ we define $A^T_{\Omega}$ as
\begin{equation*}
A^T_{\Omega}=\{(x,t)\,:\,x\in\Phi_t(\Omega),\;\;\;t\in[0,T]\}.
\end{equation*}

Setting $\mu_k= \p_x^k\rho$, $k=1, \dots, m$, the equations in \eqref{br1}, \eqref{br2} and \eqref{br3} can be written as
  \begin{equation}
 \label{a4}
 \partial_t \mu_k+a(x,t)\partial_x \mu_k = b_k(x,t)\mu_k+c_k(x,t),
 \end{equation}
 with $a(x,t)$ as \eqref{a1},
 \begin{equation*}
 b_k(x,t)\equiv (k+1)(J^{-2}-I)(\rho^2)-2\rho^2\in C([0,T] :C^1_b(\R)),
 \end{equation*}
 (see \eqref{a2}) and
 \begin{equation}
 \label{cc}
 \begin{aligned}
 c_k(x,t)&\equiv a^k_{k-1}\partial_x(J^{-2}-I)(\rho^2)\partial_x^{k-1}\rho\\
 &\;\;\;+a^k_{k-2}\partial^2_x(J^{-2}-I)(\rho^2)\partial_x^{k-2}\rho\\
&\;\;\;+\dots+ a^k_{1}\partial^{k-1}_x(J^{-2}-I)(\rho^2)\partial_x\rho\\
&\;\;\;+\partial^k_xJ^{-2}(\rho^2)\rho-\big[\partial_x^k(\rho^2)\rho-2\rho^2\partial^k_x\rho].
\end{aligned}
\end{equation}

Thus, for $k=1$,  if $\rho_0\big|_{\Omega}\in C^1$ for some open set $\Omega\subseteq \R$, since
\begin{equation*}
c_1(x,t)=\partial_xJ^{-2}(\rho^2)\rho\in C([0,T] :C^1_b(\R)),
\end{equation*}
then using the equation \eqref{a4} with $k=1$, it follows that 
\begin{equation*}
\mu_1\big|_{A^T_{\Omega}}=\partial_x\rho(\cdot,\cdot)\big|_{A^T_{\Omega}}\in C.
\end{equation*}

If $\rho_0\big|_{\Omega}\in C^2$, by combining the previous case $k=1$ and the fact that 
$\rho\in C([0,T] :H^1(\R))$ it follows that
\begin{equation}
\label{ccc}
\begin{aligned}
c_2(x,t)&=3\partial_x(J^{-2}-I)(\rho^2)\partial_x\rho+\partial^2_xJ^{-2}(\rho^2)\rho\\
&\;\;\;-\big[\partial_x^2(\rho^2)\rho-2\rho^2\partial^2_x\rho]
\end{aligned}
\end{equation}
satisfies that
\begin{equation*}
c_2(\cdot,\cdot)\big|_{A^T_{\Omega}}\in C,
\end{equation*}
then using the equation \eqref{a4} with $k=2$, one concludes that 
\begin{equation*}
\mu_2\big|_{A^T_{\Omega}}=   \partial_x^2\rho(\cdot,\cdot)\big|_{A^T_{\Omega}}\in C.
\end{equation*}

 For the general case $k\in \Z^+$, this iterative argument will yield the result  if assuming that 
 \begin{equation*}
 \partial_x^j\rho(\cdot,\cdot)\big|_{A^T_{\Omega}}\in C,\;\;\;\;\;\;j=1,2,...,k,
 \end{equation*}
 one can show that 
 \begin{equation*}
c_{k+1}(\cdot,\cdot)\big|_{A^T_{\Omega}}\in C.
\end{equation*}
 But this follows directly by the explicit form of $c_{k+1}(\cdot,\cdot)$ in \eqref{cc}.

\vspace{5mm}
\noindent\underline{\bf Acknowledgements.} 
\vspace{3mm}

The first author was partially supported by CNPq and FAPERJ/Brazil. Part of this work was completed while the second author was visiting IMPA at Rio de Janeiro whose hospitality he would like to acknowledge.


\end{document}